\documentclass[oneside,11pt,reqno,a4paper]{amsart}
\usepackage[T1]{fontenc}
\usepackage[utf8]{inputenc}
\usepackage[top=1cm,bottom=2cm,right=2cm,left=2cm,includehead]{geometry}
\usepackage[usenames,dvipsnames,svgnames]{xcolor}
\usepackage[numbers,sort&compress]{natbib}
\usepackage{doi}
\usepackage{graphicx}
\usepackage{autonum} 
\usepackage{hyperref} 
\usepackage{enumitem}
\usepackage{verbatim}

\usepackage{amsmath}
\usepackage{amssymb}
\usepackage{doi}
\usepackage[numbers,sort&compress]{natbib}
\usepackage{mathtools}
\usepackage{floatrow}
\usepackage[font={small}]{caption}
\usepackage[font={small}]{subcaption}

\newcommand{\theauthors}{Santiago Badia$^{\text{a,b}}$,
                         Manuel A. Caicedo$^{\text{a}}$,
                         Alberto F. Mart\'{\i}n$^{\text{b}}$ {\footnote{Corresponding author. \\
                         {\em E-mail addresses}: \texttt{santiago.badia@monash.edu},  
                         \texttt{mcaicedo@cimne.upc.edu},
                         \texttt{alberto.martin@monash.edu},
                         \texttt{principe@cimne.upc.edu}}},
                         and Javier Principe$^{\text{a,c}}$                         
                        
}
\newcommand{\theauths}{Manuel Caicedo} 

\newcommand{\theaffiliations}{
	\textsuperscript{a} CIMNE – Centre Internacional de M\`etodes Num\`erics en
	Enginyeria, \\ Esteve Terradas 5, 08860 Castelldefels, Spain. \\ [0.5em]
    \textsuperscript{b} School of Mathematics, Monash University, Clayton,
	Victoria, 3800, Australia. \\ [0.5em]
	\textsuperscript{c} Universitat Polit\`ecnica de Catalunya, Campus Diagonal Bes\`os, Campus Diagonal Bes\`os,  Av. Eduard Maristany 16, Edifici A (EEBE), 08019, Barcelona, Spain \\ [0.5em]
}

\newcommand{\thethanks}{Financial support from the European Commission under
the FET-HPC ExaQUte project (Grant agreement ID: 800898) within the Horizon
2020 Framework Programme is gratefully acknowledged. This work has been partially funded by the projects RTI2018-096898-B-I00 and ERC2018-092843 from the ``FEDER/Ministerio de Ciencia e Innovación –
Agencia Estatal de Investigación''. The authors thankfully acknowledge the computer
resources at Marenostrum-IV and the technical support provided by the Barcelona
Supercomputing Center (RES-ActivityID: IM-2019-3-0008, IM-2020-1-0002, IM-2020-2-0003). This work was supported by computational resources provided by the Australian Government
through NCI under the National Computational Merit Allocation Scheme. Financial support to CIMNE via the CERCA Programme /
Generalitat de Catalunya is also acknowledged.}
\usepackage{fancyhdr}

\usepackage{newtxtext}
\usepackage{newtxmath}

\usepackage{acronym}
\usepackage[most]{tcolorbox}

\usepackage{booktabs}

\usepackage{tikz}
\definecolor{myellow}{RGB}{255,230,128}
\definecolor{gray20}{RGB}{204,204,204}
\definecolor{mygray}{RGB}{204,204,204}
\definecolor{mygreen}{RGB}{138,203,95}
\definecolor{myblue}{RGB}{77,151,214}

\definecolor{lstgrey}{rgb}{0.95,0.95,0.95}
\usepackage{listings}
\usepackage[textsize=small]{todonotes}

\usepackage[vlined]{algorithm2e}
\usepackage{xcolor}

\SetAlFnt{\small}
\SetArgSty{textnormal}

\SetAlCapSty{xAlCapSty}

\SetCommentSty{xCommentSty}

\SetNlSty{mynlfont}{}{} 

\LinesNumbered

\SetSideCommentRight

\DontPrintSemicolon

\RestyleAlgo{algoruled}

\SetInd{0.5em}{0.5em}

\acrodef{dof}[DOF]{Degree Of Freedom}
\acrodefplural{dof}[DOFs]{Degrees Of Freedom}
\acrodef{vef}[VEF]{Vertex, Edge, and Face}
\acrodefplural{vef}[VEFs]{Vertices, Edges, and Faces}
\acrodef{cg}[DG]{Continuous Galerkin}
\acrodef{dg}[CG]{Discontinuous Galerkin}
\acrodef{mpi}[MPI]{Message Passing Interface}

\acrodef{fe}[FE]{Finite Element}
\acrodefplural{fe}[FEs]{Finite Elements}
\acrodef{pde}[PDE]{Partial Differential Equation}
\acrodefplural{pde}[PDEs]{Partial Differential Equations}
\acrodef{amr}[AMR]{Adaptive Mesh Refinement and coarsening}
\acrodef{amg}[AMG]{Algebraic MultiGrid}
\acrodef{hpc}[HPC]{High Performance Computing}
\acrodef{oo}[OO]{Object-Oriented}
\acrodef{bddc}[BDDC]{Balancing Domain Decomposition by Constraints}
\acrodef{cse}[CSE]{Computational Science \& Engineering}
\acrodefplural{ilu}[ILUs]{Incomplete LU decompositions}

\acrodef{hav}[HAV]{Heat Affected Volume}
\acrodef{fcm}[FCM]{Finite Cell Method}
\acrodef{ccm}[CCM]{Cut Cell Method}
\acrodef{csg}[CSG]{Constructive Solid Geometry}
\acrodef{cad}[CAD]{Computer Aided Design}
\acrodef{cae}[CAE]{Computer Aided Engineering}
\acrodef{stl}[STL]{STereoLithography}

\acrodef{sc}[SC]{simple cubic}
\acrodef{fcc}[FCC]{face-centered cubic}
\acrodef{bcc}[BCC]{body-centered cubic}

\acrodef{kee}[KEE]{Kelly Error Estimator}
\acrodef{am}[AM]{Additive Manufacturing}
\acrodef{sfc}[SFC]{Space-Filling curve}
\acrodefplural{sfc}[SFCs]{Space-Filling curves}
\acrodef{stl}[STL]{STereoLithography}
\acrodef{QoI}[QoI]{quantity of interest}
\acrodef{bvp}[BVP]{Boundary Value Problem}
\acrodef{ls}[LS]{Level-Set}
\acrodef{dof}[DOF]{Degree Of Freedom}
\acrodefplural{dof}[DOFs]{Degrees Of Freedom}
\acrodef{vef}[VEF]{Vertex, Edge, and Face}
\acrodefplural{vef}[VEFs]{Vertices, Edges, and Faces}
\acrodef{cg}[CG]{Continuous Galerkin}
\acrodef{dg}[DG]{Discontinuous Galerkin}
\acrodef{fem}[FEM]{Finite Element Method}
\acrodef{agfe}[agFE]{Aggregated Finite Element}
\acrodef{agfem}[AgFEM]{Aggregated Finite Element Method}
\acrodef{hagfem}[$h$-AgFEM]{$h$-adaptive Aggregated Finite Element Method}
\acrodef{cgm}[CG]{Conjugate Gradient}
\acrodef{mgre}[MGRE]{Maximum Global Relative Error}
\acrodef{gre}[GRE]{Global Relative Error}
\acrodef{nac}[NAC]{Number of Adaptivity Cycles}
\acrodef{cpu}[CPU]{Central Processing Unit}
\acrodefplural{cpu}[CPUs]{Central Processing Units}

\def\FEMPAR{{\texttt{FEMPAR}}}

\def\p4est{{\texttt{p4est}}}
\def\t8code{{\texttt{t8code}}}

\def\petsc{{\texttt{PETSc}}}
\def\SNES{{\texttt{SNES}}}

\def\KSP{{\texttt{KSP}}}

\def\gamg{{\texttt{GAMG}}}

\def\qhull{{\texttt{Qhull}}}
\def\mgre{{\texttt{MGRE}}}
\def\mgre{\eta_{G}^{\text{max}}}
\def\gre{\eta_{G}}

\def\mixedup{mixed (u/p)}

\def\NAS{\texttt{num\_amr\_steps}}
\def\CLS{\texttt{current\_load\_step}}
\def\CAS{\texttt{current\_amr\_step}}
\def\NLS{\texttt{num\_load\_steps}}

\def\ASF{\texttt{amr\_step\_freq}}

\def\VoidStr{hollow}

\DeclareTextFontCommand{\vbt}{\ttfamily\hyphenchar\font=45\relax}
\def\grad{{\boldsymbol{\nabla}}}

\def\Eq#1{(\ref{eq:#1})}
\def\Fig#1{Fig.~\ref{fig:#1}}
\def\Alg#1{Alg.~\ref{alg:#1}}

\def\Sec#1{Sect.~\ref{sec:#1}}

\def\Tab#1{Table~\ref{tab:#1}}

\def\A{{\mathcal{A}}}
\def\B{{\mathcal{B}}}

\def\L{{\mathcal{L}}}

\newcommand{\added}[1]{{\leavevmode #1}}

\def\bal{\begin{align}}
\def\eal{\end{align}}

\def\ph{p_h}

\def\intvars{\alpha}
\def\intvarsgen{\bs \mu}
\def\intvarsh{\intvars_h}

\def\intvarstimeh#1{\intvars^{#1}_h}
\def\intdom#1{\int_{\Dom}{#1}\,d\Dom}

\def\intgammad#1{\int_{\Gamma_{\rm D}}{#1}\,d\Gamma_{\rm D}}

\def\intgammat#1{\int_{\Gamma_{\rm N}}{#1}\,d\Gamma_{\rm N}}

\def\IdenTwo{\textbf{I}}                       

\def\SSCS{\bs{\sigma}}                         
\def\SSCSDEV{\bt{s}}                           
\def\SSPE{\bs{\varepsilon^{p}}}                
\def\SSVTE{\varepsilon_{v}}                    

\def\dev#1{\textbf{dev}(#1)}
\def\grad#1{\boldsymbol{\nabla} \, #1}

\def\trace#1{\mathbf{tr}(#1)}
\def\div#1{\nabla\cdot#1}

\def\normregl2#1#2{{\| #1 \|}_{{#2}}}
\def\normregh1#1#2{{\| #1 \|}_{1,#2}}
\def\seminormregh1#1#2{{| #1 |}_{1,#2}}

\def\zero{\mathbf{0}}
\def\NitshFact{\beta}

\def\resvector{\boldsymbol{\mathcal{R}}}

\def\cell{T}

\def\pel{\mathcal{\cell}}

\def\triang{\pel_h}

\newcommand{\cdef}{}

\def\Dom{{\cdef\Omega}}

\def\DomArt{{\cdef\Omega^{\text{art}}}}

\def\R{\mathbb{R}}

\def\x{\boldsymbol{x}}

\def\bs#1{\boldsymbol{#1}}
\def\bt#1{\textbf{#1}}

\def\u{\bs{u}}
\def\uh{\bs{u}_h}

\def\duh{\delta \uh}

\def\vh{\bs{v}_h}

\def\ph{p_h}
\def\qh{q_h}

\def\n{{\bs{n}}}

\def\f{\bs{f}}
\def\t{\bs{t}}

\def\x{\bs{x}}

\def\bou{{\Gamma}}
\def\bouD{{\bou_{\rm D}}}
\def\bouN{{\bou_{\rm N}}}

\def\L2#1{L^2(#1)}

\def\V{{\mathcal{V}}}

\def\Q{{\mathcal{Q}}}

\def\trialsp{\mathbf{V}}

\def\trialsph{\trialsp_h}

\def\ag{{\mathrm{ag}}}
\def\agd{{\mathrm{-,ag}}}
\def\act{{\mathrm{act}}}
\def\cut{{\mathrm{cut}}}
\def\integ{{\mathrm{int}}}

\def\std{{\mathrm{std}}}

\newcommand{\thetitle}{\uppercase{A robust and scalable unfitted adaptive finite element framework for nonlinear solid mechanics}}
\hypersetup{
    bookmarks=true,
    breaklinks=true,
    bookmarksopen=true,
    pdftitle={\thetitle},  
    pdfauthor={\theauths}, 
    colorlinks=true,       
    linkcolor=black,       
    citecolor=blue,        
    filecolor=black,       
    urlcolor=blue          
}

\begin{document}

\thispagestyle{empty}

\renewcommand*{\thefootnote}{\fnsymbol{footnote}}

\begin{center}
{ \bf {\thetitle}}

\vspace*{1em}

\theauthors

\vspace*{1em}

\theaffiliations

\end{center}

\setcounter{footnote}{0}
\renewcommand*{\thefootnote}{\arabic{footnote}}

\begin{center}

{\bf Abstract}

\vspace*{1em}

\begin{minipage}{0.9\textwidth}
\begin{small}
  In this work, we bridge standard \ac{amr} on scalable octree background meshes and robust unfitted \ac{fe} formulations for the \emph{automatic} and efficient solution of large-scale nonlinear solid mechanics problems posed on complex geometries, as an alternative to standard body-fitted formulations, unstructured mesh generation and graph partitioning strategies. We pay special attention to those aspects requiring a specialized treatment in the extension of the unfitted \ac{hagfem} on parallel tree-based adaptive meshes, recently developed for linear scalar elliptic problems, to handle nonlinear problems in solid mechanics. 
  In order to accurately and efficiently capture localized phenomena that frequently occur in nonlinear solid mechanics problems, we perform pseudo time-stepping in combination with $h$-adaptive dynamic mesh refinement and re-balancing driven by a-posteriori error estimators. The method is implemented considering both irreducible and \mixedup{} formulations and thus it is able to robustly face problems involving incompressible materials. In the numerical experiments, both formulations are used to model the inelastic behavior of a wide range of compressible and incompressible materials. First, a selected set of benchmarks are reproduced as a verification step. Second, a set of experiments is presented with problems involving complex geometries. Among them, we model a cantilever beam problem with spherical hollows distributed in a \ac{sc} array. This test involves a discrete  domain with up to 11.7M \acp{dof} solved in less than two hours on 3072 cores of a parallel supercomputer.
\end{small}

\end{minipage}
\end{center}

\vspace*{1em}

\noindent{\bf Keywords:} Nonlinear Solid Mechanics~$\cdot$~Adaptive Mesh Refinement~$\cdot$~Unfitted finite elements~$\cdot$~Embedded boundary methods~$\cdot$~Tree-based meshes~$\cdot$~Parallel computing.

\section{Introduction} \label{sec:introduction}

Meeting the \ac{cae} demands of many industrially-relevant settings
nowadays involves the solution of ever increasing computationally intensive problems.
Problems posed on complex geometries, which 
may evolve over time, are routinely encountered. To further increase the challenge, 
the physical phenomena subject to analysis typically exhibit 
localized features, so that the use of adaptive meshing techniques becomes paramount towards achieving an optimal accuracy versus computational efficiency balance.

A paradigmatic application problem is \ac{am}. 
\ac{am} enables complex designs with highly desirable properties largely unachievable by 
conventional manufacturing processes. For instance, \ac{am} is well-suited for producing
parts with complex mesoscale lattice structures.
The shape and interconnection pattern of unit cells hold a large influence over the mechanical properties of such structures (e.g. stress-strain response) \cite{Echeta2020}. The design of this microstructure to obtain the desired mechanical properties can be made using topology optimization \cite{Plocher2019}.
Another example is the thermo-mechanical simulation of metal \ac{am} processes, which requires that one designs a virtual mechanism that generates the growing geometry of the part being
manufactured following the real scanning path of the machine. Besides, very fine resolution close to the heated moving head is required in order to accurately capture the high thermal gradients  inherent to \ac{am} \cite{Neiva2019}, which induce distortion and residual stresses.

All this complexity limits the applicability of \ac{cae} tools traditionally used in industrial settings. These tools are almost invariably based on \acp{fe} on conforming, unstructured, body-fitted meshes. 
Their generation
for complex geometries is a challenging task, and in many cases requires human intervention. This becomes impractical when considering an external optimization loop or growing-in-time geometries, as these would require the generation of a mesh at every iteration of the process. To make things worse, parallel unstructured mesh generation and partitioning (e.g. via graph partitioning) scales poorly on parallel computers.

The approach that we advocate to effectively manage all this complexity relies on the so-called embedded (a.k.a. unfitted) \ac{fe} methods. 
The differential equation at hand is discretized by embedding the computational domain in an easy-to-generate background mesh that does not necessarily conform to its geometrical boundary, thus drastically reducing the geometrical constraints imposed on the meshes to be used for discretization.
The geometry is still provided explicitly in terms of a boundary representation (e.g. a \ac{stl} mesh)  or implicitly via the zero isosurface of a level-set function. 
Essentially, by intersecting the surface and background meshes, one generates a  sub-mesh of each cut cell that conforms to the boundary, thus generating a discretization of the domain (such discretization is only used for integration purposes, though).

Embedded methods can be used with a variety of background mesh types.
The approach herein particularly leverages octree-based meshes but the \ac{fe} schemes can readily be used for simplicial meshes. Octrees are recursively structured hexahedral grids which have multi-resolution capabilities. This is achieved by means of a recursive approach in which a mesh with a very coarse resolution (in the limit, a single cube that embeds the entire domain) is recursively refined step-by-step, until all mesh cells fulfill suitably-defined (geometrical and/or numerical) error criteria. As the terminal cells in the resulting tree-like hierarchy might be at different refinement levels, the resulting meshes are non-conforming, i.e. they have the so-called hanging \acp{vef} at the interface of neighboring cells with different refinement levels. Such relaxation of conformity becomes crucial for high parallel scalability. Besides, octree-based meshes, endowed with \acp{sfc}, enable the development of petascale-capable \ac{amr} \ac{fe} simulation pipelines,
while efficiently addressing load unbalance
caused by localization via dynamic load-balancing in the course of the simulation~\cite{Bangerth2012,BadiaMartin2019}.

On the downside, cut cells pose serious drawbacks that have reduced the applicability of embedded methods. The most concerning is that they lead to 
ill-conditioned discretizations in general. Cut cells with a small portion in the interior have a dramatic impact on the condition number of the linear system. Recently, different approaches have been considered to address this issue. A family of approaches  add stabilization terms to the discrete problem, using some kind of artificial viscosity method to make the problem well-posed; see, e.g. the CutFEM method \cite{BurmanClaus2015} or the \ac{fcm} \cite{Schillinger2015}.
More recently, the \ac{agfem} method was developed in \cite{BadiaEliptic2018,BadiaStokes2018}.
This method builds a kind of $\mathcal{C}^{0}$ Lagrangian \ac{fe} spaces that can be defined on general agglomerated meshes with arbitrary shapes. Numerical stability is achieved by eliminating the \acp{dof} laying at the exterior of the domain via suitably-defined linear algebraic constraints (in terms of a discrete extension operator), letting one stick to a Galerkin discretization. Besides, \ac{agfem} can be remarkably combined with octree-based $h$-adaptive background meshes, while being very amenable to parallelization on large-scale parallel machines~\cite{BadiaMartin2021a,Verdugo2019}. The \ac{hagfem} method, developed in \cite{BadiaMartin2021a,NeivaBadia2021}, 
combines  \ac{agfem} with parallel \ac{amr} implemented on distributed-memory platforms, and 
paves the road to functional and geometrical error-driven dynamic mesh adaptation 
with the \ac{fe} method in large-scale, industrially-relevant scenarios.

A common requirement in solid mechanics is the resolution of localized phenomena as, e.g.  in strain localization and fracture problems.
To be able to efficiently tackle these problems one has to use 
error estimators aiming at detecting the cells where the localized phenomena occurs (see, e.g. \cite{Johnson1992,Rannacher1999}). Besides, parallel \ac{amr} and dynamic load-balancing becomes necessary to efficiently address mesh densification on localized areas. In \cite{Frohne2016}, these techniques are used in
order to address problems in contact mechanics with elasto-plastic solids on parallel distributed-memory computers, while in \cite{Ghorashi2017}, a goal-oriented error estimator is developed for elasto-plasticity problems. The formulations presented in these papers are tailored to body-fitted meshes, though. 
On the other hand, unfitted formulations available in the literature are mainly restricted to linear or nonlinear elastic materials using non-adaptive meshes and serial implementations. An immersed \ac{fe} method for interface problems
was presented in \cite{Ruberg2016a} and \cite{Ruberg2016b}, using \textit{signed 
distance functions} to detect the relative position of the cell with respect to the surface 
describing the solid. The final geometry is then obtained by using \ac{csg}. 
A consistent extension of the \ac{fcm} to handle structural mechanics
problems was proposed in \cite{Schillinger2013}. The method is applied to 2D and 3D elastic
problems reproducing complex geometries imported via \ac{cad} and image-based geometric models.
A similar approach is followed in \cite{Duster2008} to reproduce 3D physical domains taken 
from human bone biopsies. The \ac{ccm} method was also extended to solve linear elasticity problems
in \cite{HansboLarson2017}. In contrast to previous works, the authors focus their attention on the
usage of parametric representations to capture the domain boundary in complex geometries.
They also extend this methodology to interface problems, in order to study the behavior of
linear elastic composite materials,  considering the embedding of thin elastic structures
such as membranes and plates.

To the best of our knowledge, and despite of this active scientific 
progress in embedded modeling, very little effort has been devoted to the design of $h$-adaptive embedded modeling methods for large-scale problems in nonlinear solid mechanics.
To fill this gap, we study the suitability of the \ac{hagfem} method~\cite{BadiaMartin2021a} on parallel scalable octree background meshes, recently developed for linear elliptic problems, to efficiently address large-scale nonlinear solid mechanics problems exhibiting localized features, posed over complex geometries. To this end, we bridge {\em for the first time in the literature} standard \ac{amr} on scalable octree (background) meshes and novel unfitted \ac{fe} methods in nonlinear solid mechanics. 
In turns out that in the extension \ac{hagfem} method~\cite{BadiaMartin2021a} to this kind of problems there are several aspects that require a somewhat specialized treatment. We restrict ourselves to problems involving nonlinear materials whose behavior depends on the strain history, i.e. \added{the constitutive models explored are based on history variables}. This sort of problems are used as a demonstrator, but
we stress that the techniques presented  can also be applied to, e.g. problems with geometric nonlinearities. In general, and even more importantly under the presence of \ac{amr}, such variables have to be tracked appropriately  in the course of a load increment simulation. To this end, we advocate for {\em considering constitutive model history variables as fields}, rather than a mere data array with values on quadrature points for numerical integration. For such functional representation, we choose a cell-wise polynomial, discontinuous Lagrangian \ac{fe} space, with the same definition of \acp{dof} for interior and cut cells, i.e. nodal values positioned at the tensor-product of Gauss quadrature points. This 
approach reduces memory demands and greatly simplifies the transfer of these variables between meshes compared to the permanent storage of history values at all quadrature points of the sub-meshes of all cut cells used for the evaluation of integrals over the embedded domain. 
Apart from this, this work relies on the following two main ingredients to achieve 
its goals:
\begin{enumerate}
\item An algorithm for the solution of strongly nonlinear problems, which composes the Newton-Raphson method with a line-search strategy with cubic backtracking. We leverage the 
suite of nonlinear and linear solvers available in the \petsc{} software package \cite{petsc-user-ref} for implementing such algorithm.
\item An algorithm to perform load increment (usually referred to as \textit{pseudo time-stepping})
  in combination with \ac{amr}. This algorithm is used to produce a locally refined background mesh
  while deformation is localizing and includes parallel dynamic load-balancing at each step.
  Up to the authors' knowledge, the current literature (see, e.g. \cite{Frohne2016,Ghorashi2017}) does not seem to pay a careful attention to this ingredient. Our experience reveals that a proper tuning of the parameters of this algorithm for each problem at hand can have a significant impact on the balance struck in practice among accuracy and computational cost.
\end{enumerate}
The discussion on these building blocks in the article focuses on the aspects that require a somewhat specialized treatment {\em in an embedded setting}.

This work is structured as follows. In \Sec{problem_statement}, we state the class of nonlinear solid mechanics problems considered in this work, namely the stress analysis of nonlinear elasto-plastic solids under small strains and displacements.  
In \Sec{agg-fe-spaces}, we introduce and motivate the spaces of functions used for spatial \ac{fe} discretization and handling of history variables, resp., and the projection operators used to 
transfer \ac{fe} functions in these spaces among two consecutive hierarchically adapted 
octree-based meshes (resulting from the application of \ac{amr}). 
In 
\Sec{numerical-solution}, we present a description of the pseudo-time 
discretization and the corresponding linearization of the nonlinear problem. 
In \Sec{mesh-adaptivity}, we describe the algorithm to bridge \textit{pseudo}-time 
integration and \ac{amr} in an embedded setting. In 
\Sec{experiments}, we present a comprehensive numerical study to validate the framework, including a set of standard validation benchmarks, and large-scale experiments involving 
complex geometries, aiming at studying the accuracy and parallel scalability 
of the framework. Finally, in \Sec{conclusions} we present some conclusions.

\section{Problem statement}\label{sec:problem_statement}

In this section we present the class of nonlinear solid mechanics problems considered herein,
which include those involving nonlinear elasto-plastic materials. 
Elasto-plasticity models can be used to solve a wide range of industrial problems.
Because they are well-known, we provide a succinct description 
of them and refer the reader to, e.g. \cite{SimoHughes1998,DeSouza2008}, 
for further details.
As usual, we use regular characters for scalar fields and bold characters for vector and tensor fields.
The two main ingredients of the model, namely the equilibrium equations and the constitutive models,
are presented in \Sec{equilibrium-equation} and \ref{sec:constitutive_model}, resp.

\subsection{Equilibrium equations}\label{sec:equilibrium-equation}

We consider both compressible and incompressible materials in this work.
The former kind of materials can be accurately modeled using an irreducible formulation, whereas the latter requires a mixed displacement-pressure formulation in which the pressure $p$ is computed separately.

Let $\Dom \subset \mathbb{R}^d$ be an open bounded domain in which the problem is posed, and $\bou$ its boundary. If we denote as $\bouD$ and $\bouN$, with $\bou = \bouD \cup \bouN$, the regions of the boundary in which we impose Dirichlet and Neumann boundary conditions, resp., the irreducible formulation consists in finding $\u : \Dom \rightarrow \R^d$, the displacement field, such that:
\begin{equation}
\begin{aligned}
\div{\SSCS} + \f       & = \zero         & & \mathrm{in} && \Dom \\
\u                         & = \overline{\u} & & \mathrm{on} && \bouD \\
\qquad \SSCS \cdot \n  & = \t            & & \mathrm{on} && \bouN \ ,
\label{eq:ss_equilibrium}
\end{aligned}
\end{equation}
given $\f : \Dom \rightarrow \R^d$ the body force per unit volume, $\overline{\u} :\bouD \rightarrow \R^d$ the displacement prescribed on $\bouD$, and $\t :\bouN \rightarrow \R^d$ the traction per unit area prescribed on $\bouN$. The unit normal pointing outwards on the boundary $\bouN$ is denoted by $\n$. The stress tensor $\SSCS$ depends nonlinearly on $\u$ as
described in \Sec{constitutive_model}.

Formulation \Eq{ss_equilibrium} fails for incompressible materials.
In the \mixedup{} formulation the volumetric part of the stress tensor 
$\SSCS$ is an additional unknown and a new equation imposing mass 
conservation is considered.
The stress tensor is decomposed as $\SSCS = p\IdenTwo + \SSCSDEV$,
where $p=\frac{1}{3}\trace{\SSCS}$ and $\SSCSDEV = \dev{\SSCS}$, denote its 
volumetric and deviatoric parts, resp. The problem can be stated as 
finding the displacement $\u : \Dom \rightarrow \R^d$, and the pressure 
$p : \Dom \rightarrow \R$ satisfying \Eq{ss_equilibrium} together with $\kappa^{-1}p-\div{\u}=0$ in $\Dom$,
where $\kappa$ is the bulk modulus.

\subsection{Constitutive model}\label{sec:constitutive_model}

The constitutive model describes the relation between stresses and strains. In 
this work, we consider the J2 Von Mises isotropic elasto-plasticity model. In any 
case, the framework is applicable to other constitutive models, e.g. damage or
coupled plasticity-damage models. As is well 
known, this constitutive model accurately describes the stress-strain response 
of a wide range of ductile materials, such as, metals and fiber-reinforced 
composites. In the context of nonlinear modeling, the relation between stresses 
and displacements can be written as
\begin{equation}
  \SSCS = \SSCS(\u,\intvarsgen), \label{eq:stress-update}
\end{equation}
that is, $\SSCS$ is a tensor function defined in terms of constitutive parameters and
the projection operator of trial stresses onto the \textit{admissible stress space}
\cite{Frohne2016,Johnson1978,Rannacher1999}
(usually implemented using the so-called \textit{return mapping} algorithm proposed 
by Simo in \cite{SimoHughes1998}) and $\intvarsgen$ denotes the set of history variables.
These variables play the role 
of tracking the nonlinearity of the constitutive model.
Particularly, in J2 plasticity, since the plastic strain is considered isochoric,
i.e. the volumetric part of the plastic deformations is zero, so the total
volumetric deformation $\SSVTE$ is directly the elastic volumetric deformation,
$\div{\u}$. The pressure $p$
is either an unknown (in the mixed formulation) or a dependent variable computed as 
$p = \kappa \div{\u}$ (in the irreducible formulation).
In the case of the mixed formulation $\SSCS$ also depends on the pressure $p$, that is,
\Eq{stress-update} becomes $\SSCS = \SSCS(\u,p,\intvarsgen)$.

The model is particularized by introducing a yield function that defines the 
limit of the elastic behavior, based on a given measure. In the case of 
J2-like models, this measure corresponds to the second deviatoric invariant 
of stresses. Concerning the set of history variables, in the isotropic case, 
this set is restricted to only one variable $\intvarsgen=\alpha$. In order to
particularize the constitutive model, we consider the following yield function \added{
\begin{equation}
  \phi = \frac{1}{2} \sqrt{\SSCSDEV : \SSCSDEV} - \sqrt{\frac{2}{3}}(\sigma_y - q(\alpha)), \label{eq:yield-function}
\end{equation}
where 
 $\sigma_y$ denotes the yield stress, and $q$ the stress-like thermodynamic  force. 
 The latter is expressed as a function of the historical variable $\alpha$ 
with an exponential saturation law as}
\begin{equation}
  q(\alpha) = -\theta H\alpha - (K_{inf}-K_0)[1-\text{exp}(-\delta\alpha)],
\end{equation}
where $\theta$ is the activation parameter to take into account the linear 
isotropic hardening, $H$ denotes the hardening modulus, and $K_{inf}$, $K_0$ and 
$\delta$ are material constants. Without loss of generality, we restrict 
ourselves to \textit{associative} elasto-plasticity models, where the yield 
function $\phi$ is taken as plastic flow rule. In consequence, the evolution 
of the plastic strain tensor $\SSPE$, and the history variables will depend on the 
choice of the yield function, i.e. \Eq{yield-function}; see \cite{SimoHughes1998} 
for further details.

\section{Spatial \ac{fe} discretization} \label{sec:agg-fe-spaces}

Unfitted \ac{fe} methods pose problems to the numerical integration and lead 
to ill-conditioned systems \cite{Burman2012,BurmanClaus2015,BadiaEliptic2018}.
To address these issues, different techniques have been proposed and our goal here is to
extend them to solve nonlinear solid mechanics problems. Specifically,
we propose the extension of the \ac{hagfem} in \cite{BadiaMartin2021a}.
We combine the conforming aggregated \ac{fe}
spaces in \cite{BadiaMartin2021a} to represent the state variables (i.e. either the displacement field or the displacement and pressure fields, depending on the formulation at hand) and discontinuous aggregated \ac{fe} spaces (see, e.g. \cite[]{muller2017high})  
to represent internal variables.

After some general aspects of the embedded boundary setup in 
\Sec{embedded-setup-cell-char}, we review the construction in \cite{BadiaMartin2021a}
in \Sec{lagr-fe-spaces} and we present the discontinuous aggregated \ac{fe} spaces
we use to represent internal variables in \Sec{aggr-disc-fe-spaces}.
Finally, in \Sec{projection-of-cm-vars} we 
describe how to use the tools provided by these spaces to define the transfer 
operators between meshes, which are required to project all variables 
involved in the problem after the mesh is adapted.

\subsection{Embedded \ac{fe} setup} \label{sec:embedded-setup-cell-char} 

\begin{figure}[htbp]
\centering
\begin{subfigure}{0.3\textwidth}
  \centering
  \includegraphics[width=0.8\textwidth]{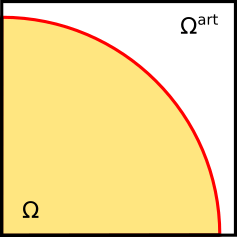}
  \caption{}
  \label{fig:setup-domain}
\end{subfigure}
\begin{subfigure}{0.3\textwidth}
  \centering
  \includegraphics[width=0.8\textwidth]{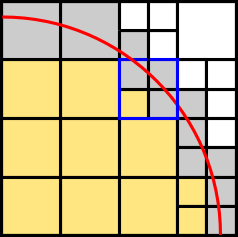}
  \caption{}
  \label{fig:setup-mesh}
\end{subfigure}
\begin{subfigure}{0.3\textwidth}
  \centering
  \includegraphics[width=0.8\textwidth]{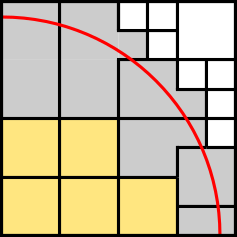}
  \caption{}
  \label{fig:setup-aggregates}
\end{subfigure}
\caption{(A) Physical domain (unfitted boundary in red) embedded into the artificial domain $\DomArt$, (B)  
         tree-based background mesh $\triang$ intersected with the physical domain in order to 
         define the \textit{well-posed} (in yellow), \textit{ill-posed} (in gray) 
         and \textit{exterior} (in white) cells, (C) Aggregated 
         mesh, \textit{well-posed} cells (in yellow) and \textit{aggregates} (in gray).}
  \label{fig:setup}
\end{figure}

Let us consider an \textit{artificial} (or \textit{background}) cuboid-like domain $\DomArt$, 
in which the \textit{physical} domain is embedded, i.e. $\Dom \subset \DomArt$ (see \Fig{setup-domain}). \added{We denote by $\triang$ a mesh for $\DomArt$.} In particular,
we consider the so-called \emph{forest-of-trees}~\cite{BadiaMartin2019} background meshes as the choice for $\triang$.
In a nutshell, these meshes are built as follows. 
First, one builds a coarse conforming mesh $\mathcal{T}_{0}$  (i.e. a quadrilateral mesh 
in 2D or hexahedral mesh in 3D) of $\DomArt$. Next, each element of this coarse mesh is the root of 
a tree-based mesh (quadtree or octree) that is generated as a result of a sequence of hierarchical refinement/coarsening steps. In our implementation we rely on
\textit{forest-of-trees} meshes generated using the \p4est{} library \cite{burstedde_p4est_2011}, which provides parallel peta-scalable mesh manipulation operations.

To describe the geometry of the \textit{physical} domain $\Dom$, let us introduce now the 
immersed boundary setting on top of the artificial domain $\DomArt$. The domain $\Dom$  (and its boundary $\partial \Dom$)  
can be represented by a level-set function or an oriented surface mesh (e.g. an \ac{stl} mesh). 
The intersection of each cell $K \in \triang$, with $\partial \Dom$ 
and $\Dom$ can be computed using a marching-cubes algorithm. We represent with $\triang^{\act}$ (resp. $\triang^{\cut}$) the set of 
active (resp., cut) cells in $\triang$ that intersect $\Dom$ (resp., $\partial \Dom$); 
$\triang^{\act} \setminus \triang^{\cut}$ is the set of interior cells in $\Dom$. 
Numerical integration over cells $K \in \triang^{\cut}$  can be carried out by computing a
simplicial decomposition of $K \cap \Dom$, e.g. using a Delaunay triangulation.
 The triangulation which results from replacing cells $K \in \triang^{\cut}$ by their decomposition is denoted by $\triang^{\integ}$. Other integration techniques that do not require simplicial sub-meshes, e.g. cubatures on general polytopes (see, e.g. \cite[]{Chin2015Dec}), can also be used.

These ingredients are enough to implement an embedded method but, as 
it was mentioned in \Sec{introduction}, a solution for the small cut cell problem is 
required. The small cut cell problem appears when the ratio between the volume of the cell
inside the domain and its total volume goes to zero. In this work, we consider as (potentially)
ill-posed cells all cut cells and well-posed cells the interior cells; see \Fig{setup-mesh}. A cell aggregation map is 
constructed to solve the small cut 
cell problem that can potentially appear in ill-posed cells; it is used to eliminate problematic \acp{dof}, as it 
will be described in \Sec{aggr-lagr-fe-spaces}. This map assigns an interior (well-posed) 
cell to any cut (ill-posed) cell, located at the 
boundary of the physical domain. In order to define this map, cell aggregates
are generated using the Algorithm 2.2 in \cite{Verdugo2019}.
Each aggregate is a connected set, composed of several ill-posed cells 
and \textit{only} one well-posed root cell $\cell$ and they form another 
partition $\triang^\ag$ (a so-called agglomerated mesh) of the domain $\DomArt$; see \Fig{setup-aggregates}.

\subsection{State variables discretization (via continuous Lagrangian \ac{fe} spaces)} \label{sec:lagr-fe-spaces}

In this section we overview two possible choices of \ac{fe} spaces for the discretization of the state variables. As required by the Galerkin \ac{fe} discretization of the problem at hand, these spaces are $H^{1}$-conforming. We use them for the  approximation of (each component of) the displacement field (and for the pressure field in the Taylor-Hood stable mixed approximation).

\subsubsection{Standard (ill-posed) continuous Lagrangian \ac{fe} spaces}\label{sec:std-lagr-fe-spaces} 

Standard {\em continuous} Lagrangian \ac{fe} spaces can be defined as:
\begin{equation}
  \V_{h}^{\std}(\triang^{\act}) \doteq \{ v \in {\mathcal{C}^0}(\Dom) \,  : \, v|_\cell \in  \V_k(\cell) \text{ for any } \cell \in \triang^{\act} \},
\end{equation}
where $\V_k(\cell)$ stands for the space of functions defined on 
$\cell \in \triang^{\act}$. Here, we consider $\V_k(\cell) \doteq \Q_k(\cell)$, the tensor product space of order $k$ univariate polynomials. Besides, we assume that all cells in $\triang^{\act}$ have local spaces $\V(\cell)$ of the same order $k$.

When using a conforming mesh, the inter-cell continuity required by the definition of $\V_{h}^{\std}(\triang^{\act})$, is implemented by using the nodal Lagrangian basis for $\Q_{k}(\cell)$ (the \acp{dof} being the corresponding nodal values) and a local-to-global \ac{dof} map. However, in general, a tree-based background mesh $\triang$  is a non-conforming mesh.
It contains the so-called hanging \acp{vef}, occurring at the 
interface of neighboring cells with different refinement levels. In this 
case, \acp{dof} lying on hanging 
\acp{vef} cannot have an arbitrary value. They must be constrained to 
guarantee trace continuity across cell interfaces. In order to simplify the computation and application of constraints, especially in a distributed-memory setting, we stick to 2:1-balanced meshes, i.e. the relation between the refinement levels of neighboring
cells is, at most, 2:1 (a precise definition can be found in  \cite{BadiaMartin2019}).

\subsubsection{Aggregated (well-posed) continuous  Lagrangian \ac{fe} spaces} \label{sec:aggr-lagr-fe-spaces}

The space $\V_h^\std$ introduced in \Sec{std-lagr-fe-spaces} is conforming, 
but leads to arbitrary ill-conditioned systems of linear algebraic equations 
due to the small cut problem mentioned before. To fix this issue, we consider 
the \ac{hagfem} recently proposed in \cite{BadiaMartin2021a}. The idea is 
to add additional constraints to $\V_{h}^{\std}$ to fix this issue, relying on 
the agglomerated mesh $\triang^{\ag}$. 

First, we consider the case of $\mathcal{C}^{0}(\Dom)$ \ac{agfe} spaces for a conforming mesh $\triang$ \cite{BadiaEliptic2018}. Let us introduce some notation. Since $\V_{h}^{\std}$ is a nodal Lagrangian \ac{fe} space, there is a one-to-one map between shape functions, nodes and \acp{dof}. The shape functions at each cell $\cell$ for $\V_{h}^{\std}$ are the standard nodal Lagrangian shape functions
$\left\{ \phi^{i}_{\cell} \right\}$, with $\V_{k}(\cell) = \mathrm{span}\left(\left\{ \phi^{i}_{\cell} \right\}_{i=1}^{n_{\Sigma}}\right)$, and $n_{\Sigma}$ the dimension of $\V_{k}(\cell)$. Each shape function $\phi^i_\cell$ is associated to a Lagrangian node $\x^i_\cell$ (defined as the coordinates of the vertices), and its corresponding \ac{dof} is $\sigma^i_\cell(v) \doteq v(\x^i_\cell)$.

In order to illustrate how $\V_{h}^{\ag}$ is built, let us start with a conforming mesh $\triang$.
We define as ill-posed the global \acp{dof} that only belongs to cut cells. We also define the aggregate that \emph{owns} this ill-posed \ac{dof} (node) among all the aggregates containing it. 
Finally, at every aggregate $\cell \in \triang^{\ag}$, we define a basis for the local aggregated
space as follows: 1) we include first the basis $\left\{ \phi^{i}_{\mathrm{root}(\cell)} \right\}_{i=1}^{n_{\Sigma}}$ for $\Q_{k}(\cell)$, i.e. the Lagrangian basis in the root cell of $\cell$,
and then we add 2) the shape functions $\left\{ \phi^{j}_{\mathrm{ill,no(\cell)}} \right\}_{i=1}^{n_{\mathrm{ill,no}}} $ in \added{ill-conditioned cells} of the aggregate that correspond 
to global \acp{dof} owned by other aggregates. 
Such space can be readily implemented by constraining the ill-posed global \acp{dof} $\left\{ \sigma_{\mathrm{ill,o}(\cell)} \right\}_{i=1}^{n_{\mathrm{ill,o}}}$ (determined by their corresponding nodes  
$\left\{ \x_{\mathrm{ill,o}(\cell)} \right\}_{i=1}^{n_{\mathrm{ill,o}}}$) in the standard space as follows:
\begin{equation}
  \sigma^{j}_{\mathrm{ill,o}(\cell)} = \sum_{i=1}^{n_{\sigma}} \phi^{i}_{\mathrm{root}(T)}(\x^{j}_{\mathrm{ill,o}(\cell)}) \sigma^{i}_{\mathrm{root}(\cell)}. 
  \label{eq:constraints}
\end{equation}
Part 1) of the aggregated local space is essential for getting optimal error estimates and eliminate the cut cell problem, while part 2) is essential to keep the $\mathcal{C}^{0}(\Dom)$ continuity. We refer to \cite{BadiaEliptic2018} for a more detailed presentation of these spaces and their numerical analysis. 

On the other hand, when $\triang$ is non-conforming, the definition of the \ac{hagfem} space involves 
two different sets of constraints, the ones related to the non-conformity of the mesh 
and the ones of cell aggregation explained above. We refer the interested reader to 
\cite{BadiaMartin2021a} for a detailed exposition of these spaces. With this 
construction, we can check that $\V_h^\ag \subset \V_h^\std$. 

\begin{figure}[h]
\centering
  \includegraphics[width=0.32\textwidth]{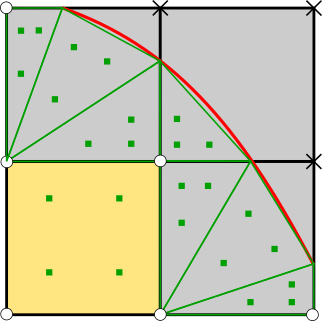}
 \caption{Illustration of the \acp{dof} of the $\V_h^\ag$ space in the patch of cells highlighted in blue in \Fig{setup-mesh} (which contains cut cells). Circles represent well-posed \acp{dof} and crosses represent ill-posed ones. Ill-posed \acp{dof} do not take arbitrary values, their values are actually defined via a suitable discrete extension operator defined in terms of linear constraints on the interior \acp{dof} of the root cell of the aggregate.  Cells $\cell \in \triang$ are shown in black and sub-meshes for integration are shown in green.  Small green squares represent numerical integration quadrature points.
 The \acp{dof} and quadrature points shown are selected to illustrate the construction. They might not reflect their actual number (and spatial position) in the numerical experiments.}
  \label{fig:agfe}
\end{figure}

The final ingredient for the
computation of \ac{fe} operators is a numerical integration in cut cells. In this work, this is done using the sub-triangulation in $\triang^{\integ}$ for each $\cell \in \triang^{\cut}$, with an appropriate choice
of the quadrature rule on simplices, see \cite{BadiaVerdugo2018} for further details.
The $\V_h^\ag$ space, together with the setup of data structures required for numerical integration in cut cells, is illustrated in \Fig{agfe}.

\subsection{History variables discretization (via discontinuous Lagrangian \ac{fe} spaces)} \label{sec:aggr-disc-fe-spaces}

\begin{figure}[htbp]
\begin{subfigure}{0.4\textwidth}
  \centering
  \includegraphics[width=0.8\textwidth]{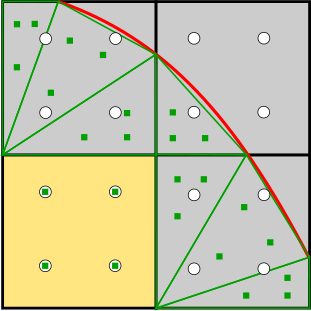}
  \caption{Discontinuous standard \ac{fe} space.}
  \label{fig:historic-spaces-without-agg}
\end{subfigure}
\begin{subfigure}{0.4\textwidth}
  \centering
  \includegraphics[width=0.8\textwidth]{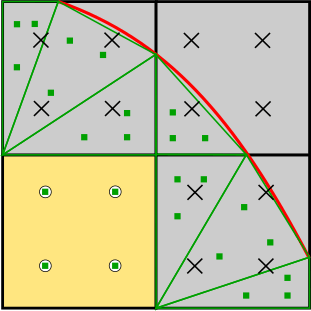}
  \caption{Discontinuous \ac{agfe} space.}
  \label{fig:historic-spaces-with-agg}
\end{subfigure}
\centering
\caption{Two possible choices for the \acp{dof} of the \ac{fe} space used to approximate the history variables. Circles represent \acp{dof} which may take arbitrary values, while crosses, \acp{dof} whose values actually depend on the interior \acp{dof} of the root cell of the aggregate via suitably-defined linear constraints.}
  \label{fig:historic-spaces}
\end{figure}

In the context of robust, $h$-adaptive, embedded \ac{fe} solvers, the representation of {\em constitutive model history variables as fields becomes essential} (for reasons made clear along the rest of the section). With a representation as a field of the history variable, we can evaluate these quantities in whatever point of the domain, rather than a mere data array with values on quadrature points for numerical integration (as it is commonly implemented in nonlinear solid mechanics codes). For such representation, we choose a cell-wise polynomial, discontinuous \ac{fe} space. In particular, the approximation of history variables is made by functions of a standard {\em discontinuous} Lagrangian \ac{fe} space on the background mesh, defined as 
\begin{equation}
  \V_{h}^{-,\std}(\triang^{\act}) \doteq \{ v \in {{L}^2}(\Dom) \,  : \, v|_\cell \in  \V_k(\cell) \text{ for any } \cell \in \triang^{\act} \}.
\end{equation}
This space is also implemented using the nodal Lagrangian basis for $\Q_{k}(\cell)$ and a local-to-global \ac{dof} map. However, as there is no required inter-cell continuity in $\V_{h}^{-,\std}(\triang^{\act})$, there is no need to glue \acp{dof} belonging to the interface of neighbouring cells. As a result, nodes can be placed anywhere within the cell (any $\x^i_\cell \in \cell$ is valid). In the case of non-conforming meshes, linear constraints on hanging \acp{vef} are not required either. Given this flexibility, we choose $\x^i_\cell$, {\em for all active cells} (i.e. both interior and cut cells), to be the location of the quadrature points used for the numerical integration in the interior cells, i.e. the tensor product of 1D Gauss quadratures; see, e.g.~\Fig{historic-spaces-without-agg}. 
Since this is the quadrature rule used for the numerical integration of the discrete operator in interior cells, in such kind of cells, the history variables are already available for numerical integration once they are computed (no extra interpolation is required). 

In cut cells it is also possible to choose a discontinuous space defined on the sub-mesh to represent history variables, that is, to consider $\V_{h}^{-,\std}(\triang^{\integ})$ as the approximating space. With this choice, the degrees of freedom in all cells (not only interior but also in cut ones) are the nodal values on the quadrature points.

However, there are a number of advantages of the representation that relies on the background mesh over this latter approach:

\begin{enumerate}
  \item First, we avoid the {\em permanent} storage of history variables in all integration points of cut cells (which become a very high number of points, specially in 3D). Instead, the value of history variables in all integration points of cut cells is interpolated locally at each cell.
  \item Second, when the background mesh is adapted (either refined and/or coarsened), a new integration mesh, and thus a new set of integration points, results from the intersection of the boundary of the domain and the adapted mesh. Thus, one needs the value of the history variables on this new set of points, while transferring functions from the original mesh to the adapted one. While this transfer operator can easily be computed in tree-based mesh refinement, it is much more involved between general unstructured meshes (as the ones in the sub-triangulation).
\end{enumerate}

  Another option available to deal with cut cells is to also use \ac{agfe} spaces with discrete extension operators for history variables. This is not required to have a well-posed problem.
  Because history variables are updated from the displacements (and eventually pressures), their \ac{dof} values are not computed solving a linear system and thus there is no actual need to classify them as well-posed or ill-posed for better conditioning. However, using an \ac{agfe} field representation we can leverage the same abstract software workflow both for the primal and internal variables, with a further reduction of computational cost in mind. Rather than computing history variables in cut cells, as it is done when $\V_{h}^{-,\std}(\triang^{\act})$ is considered (values at circles in \Fig{historic-spaces-without-agg}), they are defined by the constraints \Eq{constraints} (which give values at crosses in \Fig{historic-spaces-with-agg}).
  We have actually observed that there is a slight reduction in total computing time, as the application of the constraints in \Eq{constraints} to obtain history values in cut cells is faster than computing history variables from displacements, but the gain is a minor fraction of the total computation. We have performed a test, reported in \Sec{validation-experiments}, to verify that both options result in similar accuracy for the representation of stresses. In this case the approximating space can simply be defined as:  
\begin{equation}
  \V_{h}^{\agd}(\triang^{\ag}) \doteq \{ v \in {{L}^2}(\Dom) \,  : \, v|_\cell \in  \V_k(\cell) \text{ for any } \cell \in \triang^{\ag} \}.
\end{equation}

\subsection{Transfer operators among meshes} \label{sec:projection-of-cm-vars}

During the adaptive refinement loop, the spaces described above keep 
changing and therefore state and constitutive variables must be properly 
transferred among meshes. To this end, we use the
standard nodal interpolation and $L^2$-projection operators for refinement and coarsening, resp.,
see e.g. \cite{Ciarlet1978}.
 
However, in an embedded setting, under the presence of cut cells and aggregates, a somewhat 
specialized treatment is required. Since the refinement of a cut cell can lead to interior cells,
the nodal interpolation must also be computed for refined cut cells.
However, when refining a cut cell, 
only the interpolated \acp{dof} of the refined interior cells are used, since the (ill-posed) 
ones of the refined cut cells are defined using the constraints in  (the adapted) $\V_h^\ag$.
This applies both for \ac{fe} functions in $\V_{h}^{\ag}(\triang^{\ag})$ and $\V_{h}^{\agd}(\triang^{\ag}$). For \ac{fe} functions in $\V_{h}^{-,\std}$, we use nodal interpolation to compute  the values of the \acp{dof} of all active children cells of a refined cut cell.

\section{The nonlinear problem and its solution}\label{sec:numerical-solution}

In this section we describe the final discrete problem to be solved and 
the strategy we developed to perform this task. As usual, the nonlinear
problem is solved by gradually increasing the external loads instead of
directly looking for the final equilibrium. The load increment discretization
is described in \Sec{time-discretization} and the solution of the nonlinear
problem at each step is discussed in \Sec{linearization}.

\subsection{Pseudo-time discretization}\label{sec:time-discretization}

The final nonlinear discrete problem defined below is obtained introducing a 
pseudo-time discretization, which describes the sequential increment 
of external loads or correction of boundary conditions to reach the final 
desired ones.
Each component of the displacement field is approximated by a function 
$\uh\in \trialsph \doteq \left( \V_h^\ag \right)^d$ (the Lagrangian \ac{agfe} space of order $k$) and the internal variable is approximated by 
$\intvarsh \in \V_h^\agd$ (the discontinuous Lagrangian \ac{agfe}).

Given the history variable $\alpha_{h}^{n}$ at step $n$, the discretization of the irreducible formulation \Eq{ss_equilibrium} at load increment $n+1$ reads as follows:
find $\uh^{n+1} \in \trialsph$ and $\intvarstimeh{n+1} \in \V_h^\agd$ such that
\begin{equation} \label{eq:zero_residual_irreducible}
  \langle \resvector(\intvarstimeh{n};\uh^{n+1}) , \vh \rangle = 0,
\end{equation}
for any $\vh \in \trialsph$, where
\begin{equation} \label{eq:residual_irreducible}
  \langle \resvector(\intvarsh;\uh),\vh \rangle = \A(\intvarsh;\uh,\vh) + \B(\intvarsh;\uh,\vh) - \mathcal{F}(\intvarsh;\vh),
\end{equation}
and
\begin{align}
 \A(\intvarsh;\uh,\vh) & = \intdom{ \grad{\vh}: \SSCS(\uh,\intvarsh)}, \label{eq:residual-energy} \\ 
  \B(\intvarsh;\uh,\vh) & = \intgammad{\n \cdot  \SSCS(\vh,\intvarsh) \cdot \uh}
                - \intgammad{\n \cdot  \SSCS(\uh,\intvarsh) \cdot \vh}
                -  \intgammad{\NitshFact \vh \cdot \uh },   \label{eq:residual-boundary} \\
  \mathcal{F}(\intvarsh;\vh) & =  \intdom{\vh \cdot \f^{n+1}}  +  \intgammat{\t^{n+1} \cdot \vh} \\
  & - \intgammad{\n \cdot  \SSCS(\vh,\intvarsh) \cdot \overline{\u}^{n+1}} - \intgammad{\NitshFact \vh \cdot  \overline{\u}^{n+1}}.   \label{eq:residual-force}
\end{align}

When the mixed formulation is considered, we use the Taylor-Hood mixed \ac{fe}, but
other stable pairs can be considered; we refer to \cite{BadiaStokes2018} for the analysis
of some aggregated inf-sup stable pairs. In this case, the displacement components
space $\trialsph$  and the pressure space $\mathbf{Q}_{h}$ are 
a second order and first order aggregated Lagrangian \ac{fe} space, resp.
Then, the final discrete problem to be solved at each step is to
find $\uh^{n+1} \in \trialsph$, $p_h^{n+1} \in \mathbf{Q}_h$ and $\intvarstimeh{n+1} \in \V_h^\agd$
such that 
\begin{equation} \label{eq:zero_residual_mixed-up}
  \langle \resvector(\intvarstimeh{n};\uh^{n+1},p^{n+1}_h), (\vh,q_h) \rangle = 0 ,
\end{equation}
for any $\vh \in \trialsph$ and $q_h \in \mathbf{Q}_h$ where
\begin{equation} \label{eq:residual_mixed-up}
  \langle \resvector(\intvarsh;\uh,p_h),(\vh,\qh) \rangle = \A(\intvarsh;\uh,p_h,\vh,q_h) + \B(\intvarsh;\uh,p_h,\vh,q_h) - \mathcal{F}(\intvarsh;  \vh,q_h)),
\end{equation}
the forms $\B$ and $\mathcal{F}$ are defined
as in \Eq{residual-boundary} and \Eq{residual-force} but with $\SSCS(\uh,\ph,\intvarsh)$ and $\A$ is 
\Eq{residual-energy} plus the terms
\begin{equation}
\intdom{ \kappa^{-1}q_h p_h} - \intdom{q_h \div{\uh}}.
\end{equation}
Note that \Eq{residual-boundary} includes Nitsche's method~\cite{Schillinger2015,BurmanClaus2015,BadiaEliptic2018} terms to impose Dirichlet boundary conditions.  
$\NitshFact >0$ is a mesh-dependent parameter that needs to be large enough to ensure the positivity of the final discrete problem (after 
linearization).  For linear elliptic problems, one takes $\NitshFact =\beta_0  h_{\cell}^{-1}$ with a constant $\beta_0$ large enough. When the Dirichlet boundary is conforming to the mesh, as in the numerical experiments in \Sec{experiments}, we use a strong imposition of boundary conditions, and the Nitsche terms are switched off. The history variable $\alpha^{n+1}_{h}$ is a by-product of the evaluation of the stresses with $\uh^{n+1}$ and $\ph^{n+1}$ (see below). 

\subsection{Linearization and solution}\label{sec:linearization}

The Newton-Raphson algorithm is used to solve \Eq{zero_residual_irreducible} and \Eq{zero_residual_mixed-up}, due to its quadratic rate of asymptotic convergence. It consists of an iterative loop performed at each step in which 
the displacement is updated as
\begin{equation} \label{eq:displacement-evolution}
    \uh^{n+1,i+1} = \uh^{n+1,i} + \omega^i \duh^{i},
\end{equation}
where $i$ denotes the nonlinear solver iteration counter, $\uh^{n+1,i+1}$ 
and $\uh^{n+1,i}$ denote the value of the displacement field obtained for 
the $(i+1)-th$ and $(i)-th$ nonlinear iterations, resp., and 
$0< \omega^i \leq 1$ is a relaxation parameter determined by a standard line search
with cubic backtracking \cite{dennis1996}. 
The increment $\duh^{i}$ is obtained solving the tangent problem
\begin{equation}\label{eq:lin_sist_of_eqs}
  D\resvector(\intvarstimeh{n};\uh^{n+1,i})[\duh^{i}] = - \resvector(\intvarstimeh{n};\uh^{n+1,i}) \ ,
\end{equation}
where the Jacobian, $D\resvector$, is computed from \Eq{residual_irreducible} 
or \Eq{residual_mixed-up}, and requires differentiating $\SSCS$, 
which is implemented using the constitutive tangent tensor 
$\mathbb{C}^{\rm ep}$ see, e.g. \cite[Section 7]{DeSouza2008}.

The combination of these techniques (Newton-Raphson with line-search) is a well-known
and powerful technique for the solution of strongly nonlinear problems. Therefore, the only
point that deserves attention here is how the use of embedded methods with the choice of approximating
spaces discussed in \Sec{agg-fe-spaces} affects the evaluation of the Jacobian and the residual.
Both are required to determine the increment from \Eq{lin_sist_of_eqs}
and the residual is also required at each line search iteration.
In general, the required number of line-search iterations 
is difficult to determine a priori, as it depends on several factors, such 
as the complexity of the geometry, the loading process, and the material 
properties, among others.

The system \Eq{lin_sist_of_eqs} is assembled for the well-posed \acp{dof}. Once solved,
constrained \acp{dof} in $\duh^{i}$ are computed using \Eq{constraints}, that is, values at crosses are updated
from values at circles in \Fig{agfe}. The update in \Eq{displacement-evolution} is a functional one,
that is, performed for all \acp{dof} defining functions $\uh\in \trialsph$ (whose components are in $\V_h^\ag$).

In order to evaluate the Jacobian and/or the residual, we need to compute the stresses
(see \Eq{residual-energy} to \Eq{residual-force})
evaluating the function $\SSCS$ in \Eq{stress-update}, which includes the projection onto
the space of admissible stresses (executing the radial return algorithm).
This procedure is executed at the integration points 
(green points in \Fig{agfe} and \Fig{historic-spaces}). In interior cells, the input required for this computation is $\uh^{n+1,i}$ and $\intvarstimeh{n}$ at integration points. 
Displacements are interpolated at the integration points while history variable values at the integration points are already stored; \acp{dof} and quadrature points coincide. As a by-product of the stress projection, updated values
of the internal variable $\intvarstimeh{n+1}$ \acp{dof} are obtained. 

In cut cells, the situation is different, since history variable \acp{dof} (circles in \Fig{historic-spaces-without-agg}) are not quadrature point values (green points in \Fig{agfe} and \Fig{historic-spaces}). Thus, $\intvarstimeh{n}$ are obtained by interpolation (from circles and crosses to green points in \Fig{agfe} and \Fig{historic-spaces}, resp.). At these integration points, stresses are computed with the interpolated values, i.e. $\SSCS(\uh^{n+1,i},\alpha^{n})$.\footnote{Because the number of quadrature points in the sub-mesh generated at each cut cell can be very large (especially
  in 3D), a possible strategy to reduce this cost is to evaluate $\SSCS$ at the \ac{dof} locations of history variables
  (circles in \Fig{historic-spaces}) and then to interpolate the resulting stresses to quadrature points of the
sub-mesh (green points in \Fig{historic-spaces}). Even though the
error introduced is small, we have observed a loss of quadratic convergence in the Newton-Raphson
algorithm. Therefore, this strategy is not considered in the numerical experiments
of \Sec{experiments}.} After convergence of the nonlinear solver, the values of the history variables at the integration points (green points in \Fig{historic-spaces}) are discarded and computed at the \acp{dof} locations (circles in \Fig{historic-spaces-without-agg}) to obtain $\alpha^{n+1}$.

\section{Building hierarchically adapted meshes during incremental load-stepping}\label{sec:mesh-adaptivity}

\Alg{ilp} sketches our strategy to bridge load-stepping and dynamic \ac{amr} in {\em an embedded setting}. The algorithm is designed to 
accurately capture the inelastic behavior of the embedded solid in the course of 
a load-increment simulation, while keeping the computational demands within 
acceptable margins. To this end, the algorithm exposes a set of user-defined 
parameters that play a role in striking a balance among controlling that the global error does not blow up during the simulation and keeping computational requirements within acceptable margins.
Up to the authors' knowledge, the approaches proposed in the literature to bridge \textit{load-stepping} and \ac{amr} (see, e.g. \cite{Frohne2016,Ghorashi2017}), do not seem to pay a careful attention to this ingredient, advocating instead for simple strategies
in which the mesh is adapted either once or at most once at each step (regardless of the global error).  Our experience reveals that a proper tuning of the parameters of \Alg{ilp} for each problem at hand can have a significant impact on the balance struck in practice among accuracy and computational cost.

{
\setlength{\algomargin}{1em}
\SetFuncSty{\texbf}
\SetKwComment{comment}{}{}{}
\SetCommentSty{textrm}
\SetKwProg{Fn}{Function}{}{end}
\SetKwInOut{Input}{input}
\SetKwInOut{Output}{output}

\begin{algorithm}[htbp!]
Generate initial background mesh 
(uniformly refined octree + multiple sweeps of coarsening of exterior cells) \; \nllabel{line:initial_mesh}

\CLS $ \leftarrow$ 1 \;
\While{\CLS{} $ < $ \NLS{}}{ \nllabel{line:outer-loop}
  Update current load/displacement \; \nllabel{line:update-displ}
  Solve nonlinear problem [see \Eq{residual_irreducible} or \Eq{residual_mixed-up}] \; \nllabel{line:solve-nlp-first}
  \If{ \added{\CLS{}  modulo \ASF{}} $\equiv 0$ }{ \nllabel{line:adapt-step}
  Compute $\eta_\cell$, $\cell \in \triang^{\act}$ and $\gre$ [see \Eq{gre}] \; \nllabel{line:compute-error-est-first}
  \CAS{} $ \leftarrow$ 1 \;
  \While{($\gre$ $>$ $\mgre$) and (\CAS{} < \NAS{})}{ \nllabel{line:inner-loop}
        Retrieve state and constitutive model variables of current load step \; \nllabel{line:retrieve-prev-sol}
        
        Refine/coarsen fraction $\theta_r$/$\theta_c$ of $|\triang^{\act}|$ with the largest/smallest $\eta_T$ \nllabel{line:mark}
        
        Intersect domain with adapted background mesh and generate aggregates
        
        Transfer variables to the new mesh [see \Sec{projection-of-cm-vars}] \; \nllabel{line:adapt-mesh}
        
        Redistribute mesh and migrate variables among tasks \; \nllabel{line:redistr-tasks}

        Intersect domain with adapted background mesh and generate aggregates

        Solve nonlinear problem [see \Eq{residual_irreducible} or \Eq{residual_mixed-up}] \; \nllabel{line:solve-nlp-second}
        Compute $\eta_\cell$, $\cell \in \triang^{\act}$ and $\gre$ [see \Eq{gre}] \nllabel{line:compute-error-est-second} \; 
        \CAS{} $\leftarrow$ \CAS{} + 1 \; \nllabel{line:update-arm-counter}
  }
  } \nllabel{line:adapt-step-end}
  \CLS{} $\leftarrow$ \CLS{} + 1  \;
}
\caption{Load-stepping and parallel \ac{amr} in an embedded setting. \label{alg:ilp}}
\end{algorithm}

The inner \ac{amr} loop in \Alg{ilp} combines (a-posteriori) local and global 
error estimators. These are denoted as $\eta_\cell$, for $\cell \in \triang^{\act}$, 
and $\eta_G$, resp., in \Alg{ilp}. Local error estimators 
are used in order to decide which cells  are marked for refinement and coarsening 
among two consecutive meshes in the hierarchy; see line \ref{line:mark}. 
In particular, for that 
purpose, given user-defined refinement and coarsening fractions, denoted by 
$\theta_r$ and $\theta_c$, resp., we refine the fraction $\theta_r$ of the total
number of {\em active cells} with the largest local error estimator, and coarsen the fraction 
$\theta_c$ with the smallest local error estimators. 
We note that the generation of the initial background mesh in line \ref{line:initial_mesh} applies, to a uniformly refined octree, multiple sweeps of coarsening of exterior cells till two consecutive meshes are equivalent. This pursues to minimize the number of exterior cells of the background mesh .\footnote{While exterior cells do not carry \acp{dof}, they are handled by the underlying octree meshing engine, and thus it is convenient to reduce them to the minimum to avoid an extra overhead.}  
The experiments in \Sec{experiments} reveal that this strategy is effective in achieving large percentages of active cells versus total number of cells along the whole simulation.

In \Sec{experiments}, for simplicity,
we consider a residual-based a-posteriori error estimator grounded on the seminal works
\cite{Babuska1978,Babuska1981,Kelly1983} (although more complex local error estimators can readily be used in our framework). This kind of error estimators are well-known, so that we refer the reader to these references for a detailed definition. \added{It is worth} mentioning, however, that in an embedded setting, the implementation of this algorithm requires integration on ``cut faces'' (i.e. faces of the cut cell sub-mesh of simplices that enclose $\cell \cap \Omega$; see \Fig{agfe}), 
as these error estimators involve the evaluation of integrals (of the jump of the normal stresses) on the boundary of the cells.  On the other hand, in \Sec{experiments}, we consider the following {\em estimation} of the relative global error 
\begin{equation} \label{eq:gre}
\gre = \sqrt{\frac{ \sum_{ \cell } \eta_\cell^{2} }{ \A(\uh,\uh)} },
\end{equation}
where $\A(\uh,\uh)$ (i.e. $||\uh||_\A^2$) denotes the (squared) discrete energy norm defined in \Eq{residual-energy}. This definition of the estimated global error was proposed to drive mesh adaptivity for linear elliptic problems in \cite{Gratsch2005}.
This definition has the benefit of being dimensionless, in contrast to, e.g. using the Euclidean norm of the vector of nodal values of $\uh$ in the denominator \cite{Ghorashi2017}.

The estimated global error $\gre$, combined with a user-defined upper 
threshold, denoted as $\mgre$ in \Alg{ilp}, defines a mesh acceptability criterion that we are willing 
to fulfill at the final load, but also plays a role along the simulation by 
controlling how many adaptation steps we perform within each load step.
The requirement on $\gre$ strictly below $\mgre$ is relaxed with the so-called 
\NAS{} parameter; see line \ref{line:inner-loop}. Essentially, in each load step, a maximum of \NAS{} mesh adaptation steps are allowed. 
\added{We note that the execution of the inner loop in \Alg{ilp} depends on the user-configurable \ASF{} parameter, the frequency that controls how many load steps the user allows the simulation to run without adapting the mesh.}

\section{Numerical experiments}\label{sec:experiments}

The main purpose of this section is to showcase the high suitability of
the algorithmic framework proposed in order to efficiently deal, both in terms
of numerical accuracy and computational performance, with large-scale elasto-plasticity
problems  with localized effects posed over complex geometries. 
First, we have conducted a set of validation tests of the framework against experimental benchmarks
with results available in the literature and problems with known analytical solution.
In particular, we considered the 2D nearly-incompressible plane strain Cook's Membrane
(see, e.g. \cite{Elguedj2008}), an internally pressurized thick cylinder
(see, e.g. \cite{Ghorashi2017}  and \cite[Section 7.5.1]{DeSouza2008}) and the
stretching of a 3D perforated rectangular plate (see \cite{SimoTaylor86,DeSouza2008}). 
{\em We have been able to confirm the correctness of the algorithms at hand when applied
  to these benchmarks, for which our results match the ones in the literature accurately.}
Among these examples, \added{for conciseness}, we only present here the results of the internally pressurized thick
cylinder in \Sec{validation-experiments}. \added{In particular, this section compares the accuracy of both possibilities
for the approximation of the history variables discussed in \Sec{aggr-disc-fe-spaces}.} Before that,
we describe the experimental environment (hardware and software) and the setup common to all experiments in \Sec{exp-environment}.

\added{Moving to complex geometries, we consider in \Sec{traction} the traction of two 3D periodic lattice-type geometries representative of the ones which are enabled by \ac{am} technology; see, e.g., \cite{Tao2016}. Although there is no analytical solution at our disposal for such complex elasto-plasticity problems, the focus here is to evaluate the ability of the building
blocks of \Alg{ilp} to resolve the underlying physics while keeping the size of the discrete problem within acceptable margins. Later on, in \Sec{CBT-3D}, we consider 
the deformation of a cantilever beam with spherical hollows subject to a vertical uniformly distributed load applied in the upper face. This section also reports the results of a strong scalability study. The goal is to evaluate the ability of the 
framework to efficiently reduce computational times with increasing computational resources (cores and memory) when deployed on a parallel supercomputer.

We stress that the problems in this section have been judiciously selected so that they have 
features that potentially pose significant challenges to \ac{cad}/\ac{cae} frameworks: intricate geometries, high nonlinearity, long and thin regions with 
strong localization (weak links), and varying amount of material per unit of volume.
These features cause non-symmetric refinement patterns, 
potential load unbalance to be addressed dynamically in a parallel setting, and
nonlinearity onset on the unfitted boundary, among others. The main goal of the 
section is to reassure that the algorithmic framework 
proposed is able to promisingly deal with all this complexity.}

\subsection{Experimental environment and setup common to all experiments}\label{sec:exp-environment}

All numerical experiments are carried out on a parallel, distributed-memory environment. 
In particular, we used two different supercomputers:

\begin{itemize}
 \item {\bf Marenostrum-IV} (MN-IV) \cite{marenostrum}, hosted by the Barcelona 
 Supercomputing Centre. MN-IV is a petascale machine with 3456 nodes distributed in 
 45 racks, interconnected with the Intel OPA HPC network. Each node has 2x Intel 
 Xeon Platinum 8160 multi-core \added{cores}, with 24 \added{cores} each (i.e. 48 \added{cores} per node), and 96 GB of RAM. 
 
 \item {\bf NCI-Gadi} \cite{NCIgadi}, hosted by the Australian National 
 Computational Infrastructure Agency (NCI), is a petascale machine with 3024 nodes, 
 each containing 2x 24-core Intel Xeon Scalable \textit{Cascade Lake} \added{cores} 
 and 192 GB of RAM.
 All nodes are interconnected via Mellanox Technologies' latest generation 
 HDR InfiniBand technology.
\end{itemize}

The algorithms at hand have been implemented using the tools provided by \FEMPAR{} \cite{badiafempar2017}, 
an open source \ac{oo} Fortran200X scientific software package for the \ac{hpc} simulation 
of complex multiphysics problems governed by \acp{pde} at large scales. Among others, 
\FEMPAR{} provides scalable implementations of the algorithms in charge of generating the \ac{agfem} aggregates, and handling the linear constraints required in order to preserve conformity on non-conforming octrees and make the problem better conditioned, all in a parallel distributed framework; see ~\cite{BadiaMartin2021a} for a detailed presentation of these algorithms. We configure \FEMPAR{} so that it uses \p4est{} as its specialized forest-of-octrees meshing engine~\cite{BadiaMartin2019}. In order to solve the nonlinear
problem (see \Sec{numerical-solution}) 
we exploit the algorithms available in the \textit{Scalable Nonlinear Equations Solver}
\SNES{}, a module available in \petsc{} \cite{petsc-user-ref}. \FEMPAR{} provides
to \petsc{} all data structures required to properly handle the life cycle of the 
nonlinear solver. \FEMPAR{} \cite{BadiaFempar2020} was linked against \p4est{} v2.2 
\cite{burstedde_p4est_2011}, \petsc{} v3.12.4 \cite{petsc-user-ref} and 
\qhull{} v2015.2 \cite{qhull}, a library to compute convex hulls and Delaunay 
triangulations. All software is compiled with Intel v.18.0.5 in MN-IV, Intel v.19.1.0.166 
in {\em NCI-Gadi}. All floating-point operations are performed in IEEE double precision.

In all experiments of the section the geometry is built using a cuboid $[0,L]^d$ 
as artificial domain, 
discretized with a background octree mesh of this domain. This mesh is then intersected with level-set functions, which 
are used to describe the boundary of the embedded domain.
To keep the presentation concise, we do not include the mathematical definition of 
these functions. We refer the reader to, e.g. \cite{BurmanClaus2015,Verdugo2019},
and references therein for a detailed definition.
\Tab{params-irreducible} summarizes the main parameters and computational strategies
used in the experiments. The 
parameters of \Alg{ilp} are particularized for every test, 
and thus are presented in the section corresponding to each test. 

\begin{table}[ht!]
	\centering
	\begin{small}
		\begin{tabular}{ll}
			\toprule
			Description & Considered methods/values \\
			\midrule
			Model formulation & Irreducible formulation \vspace{0.12cm} \\
			Background mesh & Single quadtree (2D) or octree (3D) \vspace{0.12cm} \\
			\ac{fe} space state vars & First order \ac{agfe} space \vspace{0.12cm} \\
			\ac{fe} space history vars & First order discontinuous \ac{agfe} space \vspace{0.12cm} \\
			\SNES{} - nonlinear solver & Newton + cubic backtracking line search
			\vspace{0.12cm} \\		
			\SNES{} stopping criterion & $\| \mathbf{r} \|_2 < 10^{-12}$ \\
			 & $\| \mathbf{r} \|_2/ \|\mathbf{r_0} \|_2 < 10^{-12}$
			\vspace{0.12cm} \\
			\KSP{} - linear solver & Preconditioned conjugate gradient 
			\vspace{0.12cm} \\
			Parallel preconditioner & Smoothed-aggregation \gamg{}~(PETSc) \vspace{0.12cm} \\
			\KSP{} stopping criterion & $\| \mathbf{r} \|_2/ \| \mathbf{b} \|_2
			< 10^{-12}$ \vspace{0.12cm} \\
			\bottomrule
		\end{tabular}
	\end{small}
	\caption{Summary of main parameters and computational strategies.}
	\label{tab:params-irreducible}
\end{table}

\subsection{Internally pressurized thick cylinder} \label{sec:validation-experiments}

In this section we present selected results when applying the framework
to an elasto-plastic problem with known analytical solution, namely a
long metallic thick-walled cylinder subjected to internal pressure.
The main goal is to evaluate the impact 
of the \ac{fe} space used for representing the internal variables of the 
constitutive model, i.e.  $\V_{h}^{-,\std}$ versus $\V_{h}^{\agd}$,
on the accuracy of the computed stresses.
The main feature of this  test is that, once the internal pressure reaches its maximum value, 
the whole domain is in the plastic regime.

The analysis of the internally pressurized thick cylinder is carried out assuming plane strain conditions.
As a consequence, only one quarter of the cylinder, with 
the appropriate boundary conditions, is required to represent the full test 
(see \cite{Ghorashi2017} and \cite[Section 7.5.1]{DeSouza2008}); see \Fig{TTT2D-geometry}
for  the description of the geometry and boundary conditions.
The pressure $P$, prescribed in the inner surface, is incremented gradually 
until the limit load is reached. While the internal pressure is small, the entire 
cylinder will remain elastic. However, as $P$ increases, the cylinder begins to 
yield from the inner surface $r=a$. Then, the yielded region expands outwards 
forming a cylindrical plastic front. The analytical displacement solution in polar 
coordinates for both, elastic and plastic regions, as a function of $r$, can be found in 
\cite{Hill1950} and \cite{Gao2003}.

We set the inner radius $a=100 \, \text{mm}$, and outer radius $b=200 \, \text{mm}$
and the side of the artificial domain has a length $L=b$. 
In order to capture the behavior of the material, we considered a compressible 
perfect plastic material with Young's modulus $E=210 \, \text{GPa}$, Poisson's 
ratio $\nu=0.30$, and yield stress $\sigma_y=240 \, \text{MPa}$. We considered 
as total internal pressure $P=0.19\, \text{GPa}$, discretized into 11 load steps.

In \Fig{TTT2D-convergence-test} we report the $L^2$-norm relative error in 
the radial ($\sigma_{{r}{r}}$) and hoop stress ($\sigma_{{\theta}{\theta}}$) components
as a function of the number of \acp{dof} 
in the final state of the mesh. The different points in the figure were generated using 
a sequence of uniformly refined meshes, starting from a $4\times 4$ uniform mesh at the
left-most point refined six times to reach the right-most one.
\added{The error of the radial stress decays as $\mathcal{O}(h^{-1})$, i.e.
  as $\mathcal{O}(n^{-1/2})$, with $n$ the number of \acp{dof}, as it can be expected
  for an elliptic problem with a smooth solution approximated by linear functions.
  The error of the hoop stress also decays as $\mathcal{O}(n^{-1/2})$ for coarse meshes
  but the rate decreases for finer ones. Because the exact hoop stress
  is not regular \cite{Gao2003}, convergence can be expected but its rate cannot be predicted.}
In any case, the main conclusion that can be extracted from the experiment is that
both options for the \ac{fe} space representing the history variables result in
similar accuracy for the representation of stresses. As a slight reduction of computational
times was observed for $\V_{h}^{\agd}$, we stick into this \ac{fe} space for the 
rest of the experiments in the section.

\begin{figure}[h!]
\centering
\begin{subfigure}{0.45\textwidth}
  \centering
  \includegraphics[width=0.7\textwidth]{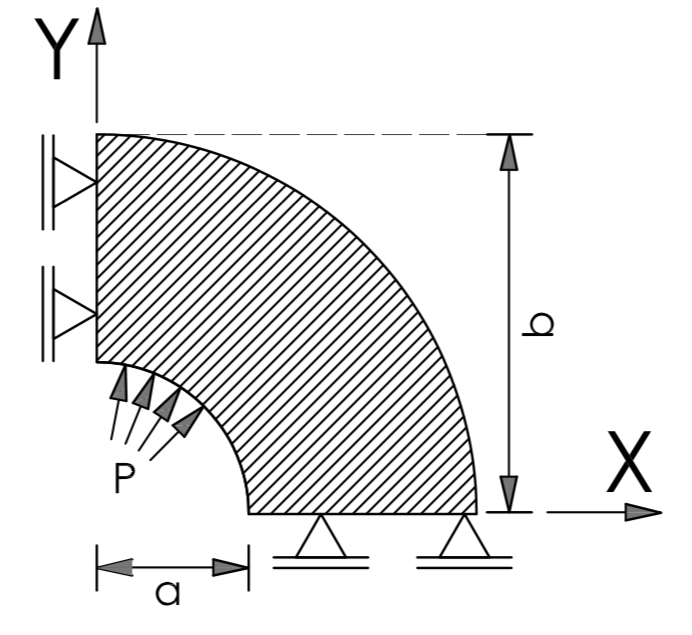}
  \caption{}
  \label{fig:TTT2D-geometry}
\end{subfigure}
\begin{subfigure}{0.45\textwidth}
  \centering
  \includegraphics[width=1.0\textwidth]{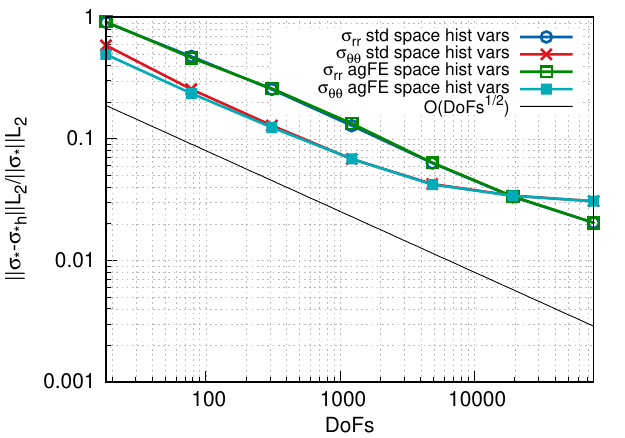}
  \caption{}
  \label{fig:TTT2D-convergence-test}
\end{subfigure}
\caption{Internally pressurized thick cylinder: (A) Geometry description and boundary conditions. (B) Convergence of radial ($\sigma_{{r}{r}}$) and hoop 
($\sigma_{{\theta}{\theta}}$) stresses. $\V_{h}^{-,\std}$ (standard space for history variables) versus $\V_{h}^{\agd}$ (\ac{agfe} space for history variables).}
  \label{fig:TTT2D-results}
\end{figure}

\subsection{Traction of complex 3D lattice-type structures}
\label{sec:traction}

\added{After validating the framework against standard benchmarks, 
we move to more intricate geometries with the aim of verifying the ability of the building
blocks of \Alg{ilp} to accurately resolve the underlying physics
while keeping computational resources (e.g. number of cells/\acp{dof}) within acceptable margins. 
To this end, we consider traction tests performed on two different 3D periodic structures
representative of the ones that are enabled by \ac{am} technology \cite{Tao2016}; 
see \Fig{Domains}.
This kind of structures are geometrically intricate, difficult to discretize with
body-fitted mesh generators.  This complexity triggers also a nontrivial physical behavior
with long and thin regions (weak links) where localization may occur. 
Besides, the onset of nonlinearity is located on the 
unfitted boundary, putting further stress on the adaptive cut cell \ac{fe} formulation.
The two geometries considered herein
differ in the amount of material per unit of volume (the ratio between material
and empty space), resulting in different refinement patterns and geometric distribution of 
plastic regions.}

\begin{figure}[h!]
  \centering
  \begin{subfigure}{0.4\textwidth}
    \centering
    \includegraphics[width=\textwidth]{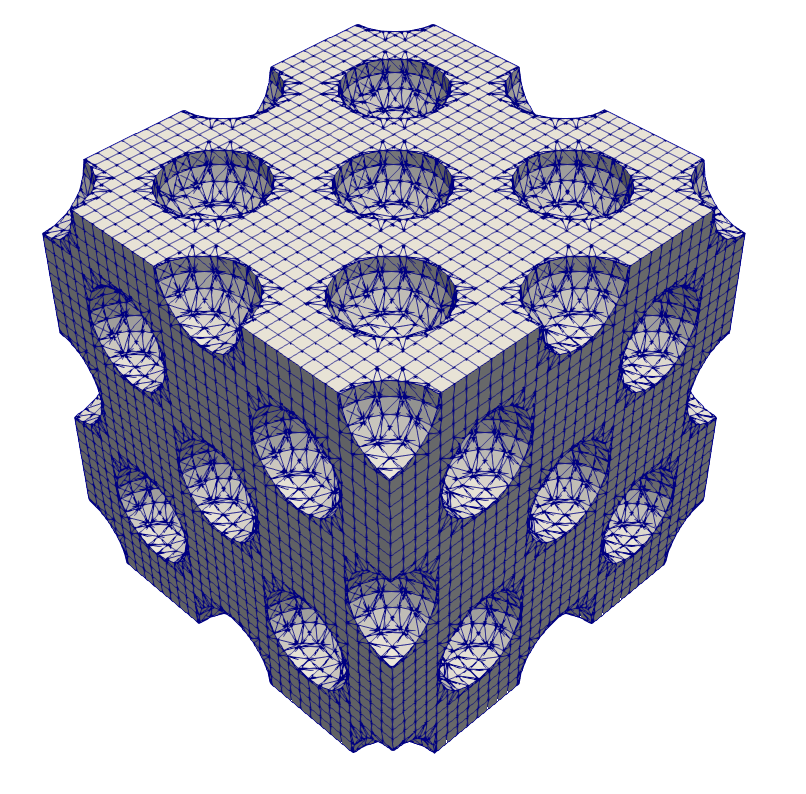}
    \caption{}
    \label{fig:DomainBulletSpheres}
  \end{subfigure}
  \begin{subfigure}{0.4\textwidth}
    \centering
    \includegraphics[width=\textwidth]{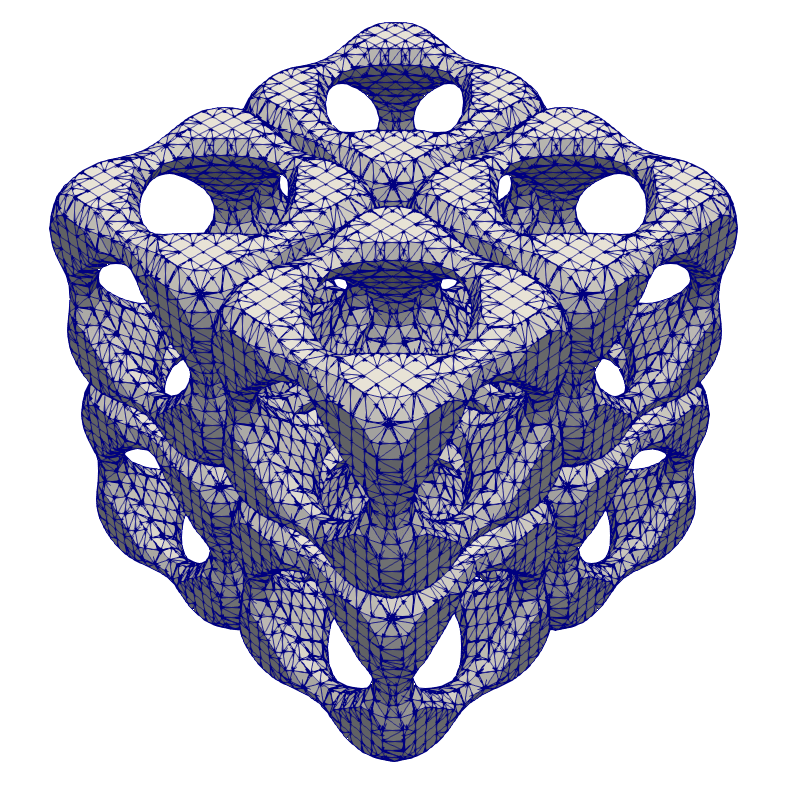}
    \caption{}
    \label{fig:DomainSwissCheese}
  \end{subfigure}
  \caption{3D geometry and initial \ac{fe} mesh ($\triang^{\integ}$) of (A) \VoidStr{} and (B) lattice structures.}
  \label{fig:Domains}
  \end{figure}

\added{
The first geometry is a sample of material that consists of homogeneous medium
with spherical hollows of radius $0.12\, \text{m}$, 
organized in a periodic \ac{fcc} lattice structure of side $L/2$.
\cite{Kittel2004}; see \Fig{DomainBulletSpheres}.
The second geometry is a lattice structure obtained by repetitive Boolean operations
(copy, translate and union) of a single structural block described by level-set
functions given in \cite{BurmanClaus2015}; see \Fig{DomainSwissCheese}.
The artificial domain is the unitary cuboid
with side length $L=1\,\text{m}$ in both cases.
This domain can be simply discretized by a single octree.
\Fig{Domains} also shows the initial $\triang^{\integ}$ meshes, i.e. 
the ones which result from the intersection of the domain and the 
initial background meshes. The latter were in turn generated by uniform refinement of
the root cell of the octree plus multiple sweeps of coarsening of exterior cells 
(see line \ref{line:initial_mesh} of \Alg{ilp}).
We note that the meshes shown include the tetrahedral subcells generated for
integration purposes (excluding exterior cells), so that 
we can clearly see the shape of the boundaries. A body-fitted discretization of these geometries is not
easy to generate but automatic with the unfitted method.}

\added{In both tests we considered the following material properties:
Young modulus $E=70 \, \text{GPa}$ and 
Poisson's ratio $\nu=0.2$. The inelastic material behavior will be 
captured by an isotropic J2 plasticity model based on the particularization 
of the functions presented in \Sec{constitutive_model}. For the hardening 
evolution law, we consider a linear hardening ($K_{inf}=K_0=0$ and 
$\delta=0$), the yield threshold and the isotropic hardening factors being
$\sigma_{y}=0.243 \, \text{GPa}$ and $H=0.2 \,  \text{GPa}$, resp.
On the other hand, boundary conditions are also similar. The structure is clamped at $x=0$ and
the rest of the faces of the structures are traction-free. An horizontal displacement 
$\u_x=0.01\,\text{m}$ is incrementally imposed in $60$ uniform load steps at $x=L$ in 
the case of the \VoidStr{} structure, while $\u_x=0.005\,\text{m}$ is
incrementally imposed in $51$ uniform load steps at the same face in the lattice structure.
We used \Alg{ilp} taking $\theta_r=0.1$ and $\theta_c=0.05$, and a constant
$\mgre=0.1$ along the whole simulation. Besides, we set \ASF=4, and \NAS{}=4.

In \Fig{Plasticity} we show the value of 
the J2 plasticity internal variable $\alpha$ once the full load has been applied to the corresponding structures. 
This is accompanied, in \Fig{Meshes}, with the associated mesh refinement pattern. 
For a better illustration, Figs. \ref{fig:MeshBulletSpheres} and \ref{fig:MeshSwissCheese} also include exterior cells. We stress, though, that these cells do not actually carry \acp{dof} in the unfitted method.}

\begin{figure}[h!]
  \centering
  \begin{subfigure}{0.4\textwidth}
    \centering
    \includegraphics[width=\textwidth]{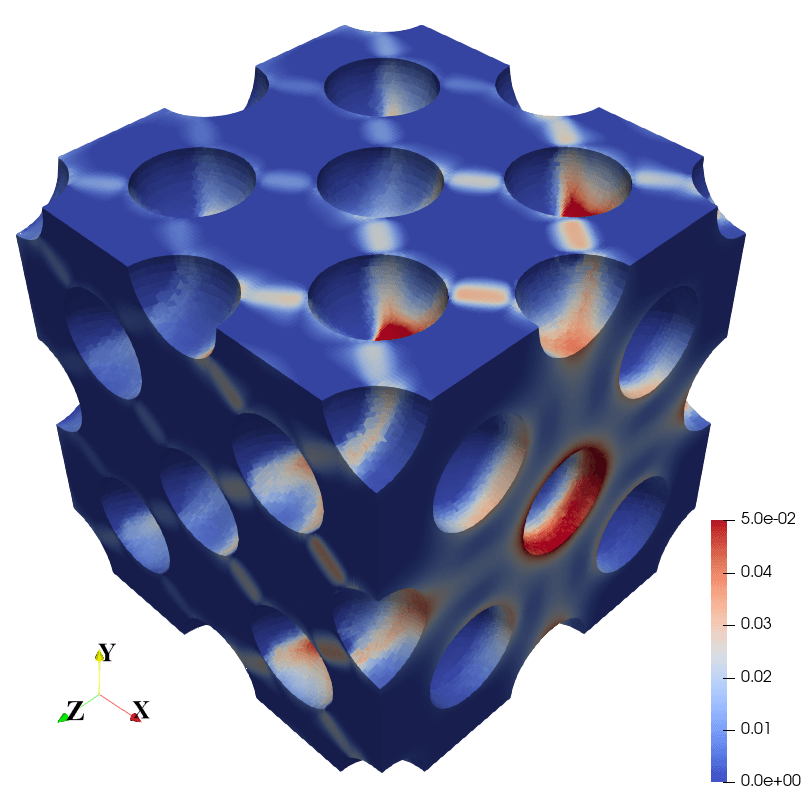}
    \caption{}
    \label{fig:PlasticityBulletSpheres}
  \end{subfigure}
  \begin{subfigure}{0.4\textwidth}
    \centering
    \includegraphics[width=\textwidth]{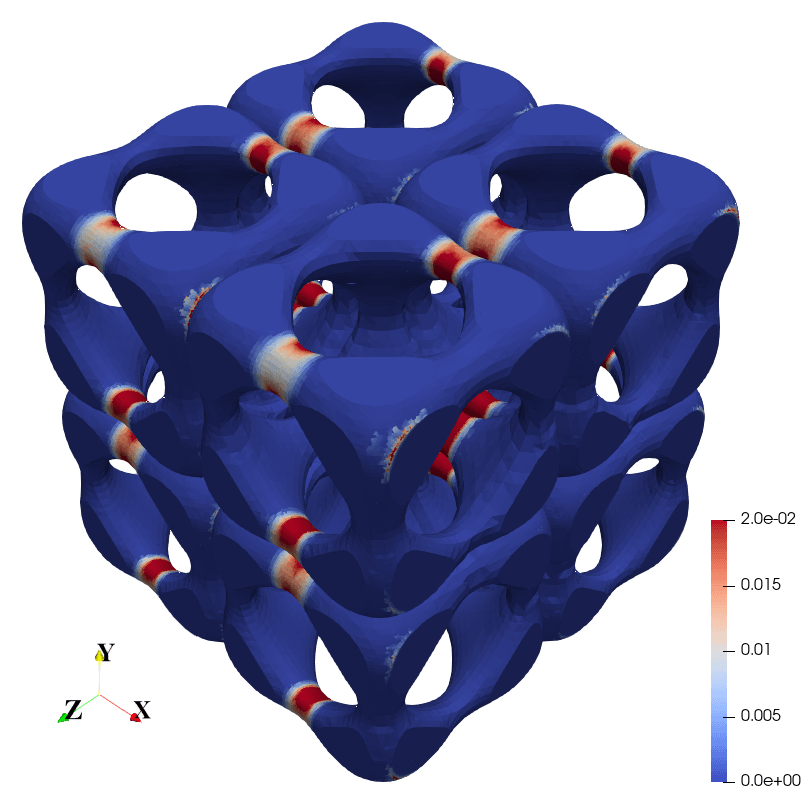}
    \caption{}
    \label{fig:PlasticitySwhissCheese}
  \end{subfigure}
    \caption{J2 Plasticity history variable $\alpha$ at the last load step  
    in (A) \VoidStr{} and (B) lattice structures.}
    \label{fig:Plasticity}
  \end{figure}
  }
  
  \begin{figure}[h!]
  \centering
  \begin{subfigure}{0.3\textwidth}
    \centering
    \includegraphics[width=\textwidth]{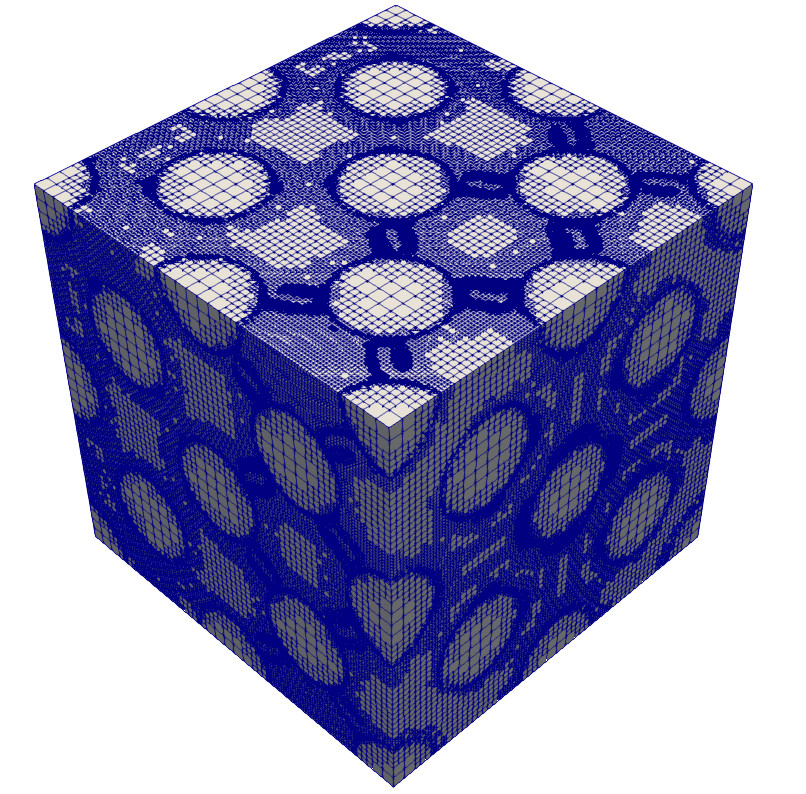}
    \caption{}
    \label{fig:MeshBulletSpheres}
  \end{subfigure}
  \begin{subfigure}{0.3\textwidth}
    \centering
    \includegraphics[width=\textwidth]{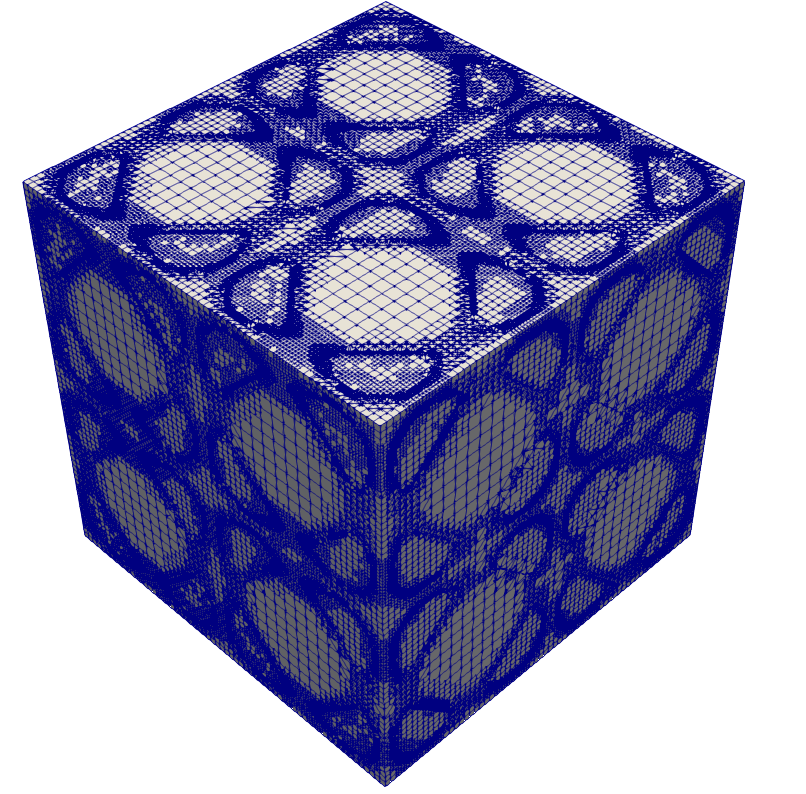}
    \caption{}
    \label{fig:MeshSwissCheese}
  \end{subfigure}
  \begin{subfigure}{0.3\textwidth}
    \centering
    \includegraphics[width=\textwidth]{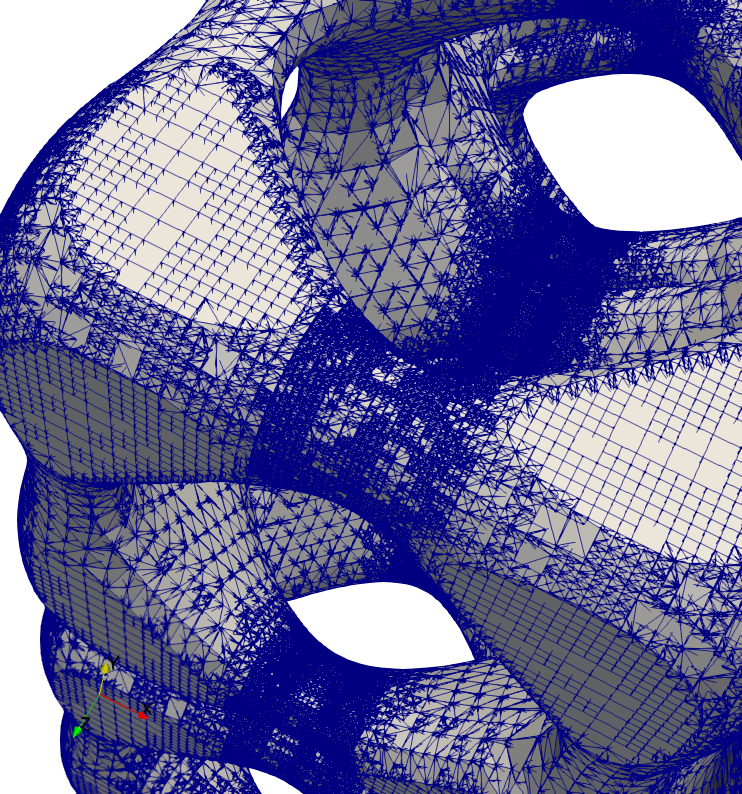}
    \caption{}
    \label{fig:MeshSwissCheeseZoom}
  \end{subfigure}
  \caption{Mesh refinement pattern in (A) \VoidStr{} and (B) lattice structures. 
  (C) Close view of the mesh of the lattice structure.}
  \label{fig:Meshes}
  \end{figure}

\begin{figure}[h!]
  \includegraphics[width=0.8\textwidth]{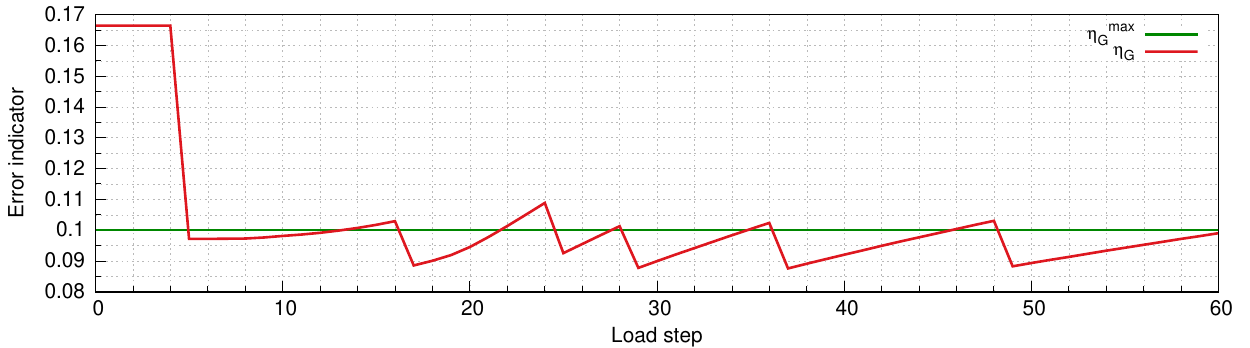}
  \caption{Evolution of the error indicator during loading (\VoidStr{} structure).}
  \label{fig:IndicatorsError}
\end{figure}

\added{Despite not having access to the exact solution, we can conclude that the 
results obtained are in a reasonable agreement with the physics underlying the problem at hand.
The variable $\alpha$ can be interpreted as a measure of the evolution of the inelastic strains of the material. Regions where the stress states are mainly controlled by shear (i.e. involving a large deviatoric component) will yield large values of $\alpha$. 
{\em \Fig{Plasticity} clearly reveals some localized regions with higher $\alpha$ values, with a particular spatial distribution which is consistent with the geometry of the structures and the applied loads.} This effect is accompanied with a higher  local mesh resolution in those areas, as can be observed in \Fig{Meshes}.  
In the \VoidStr{} structure (\Fig{PlasticityBulletSpheres}), some strain localization bands can be observed in regions between hollows. (Nevertheless the plastic regime has also spread through many other regions of the structure.) In contrast, in the lattice structure, shown in \Fig{PlasticitySwhissCheese}, the plastic region is concentrated exclusively in the thin links of the structure.

In the sequel, for conciseness, we restrict the presentation of results to those of the \VoidStr{} structure. The results obtained with the lattice structure are, in general qualitative terms, similar. In particular, in \Fig{IndicatorsError} we show the evolution of the error indicator $\gre$, while in \Tab{SummaryBullets}, we report the number of \acp{dof} and cells, and the fraction 
of active cells with respect to the total number of cells in the background mesh for each of the stages in-between mesh adaptations. The former quantity as reported in \Tab{SummaryBullets} (and in any of the tables herein)  does not include those \acp{dof} which are subject to multipoint linear constraints (i.e. hanging and ill-posed \acp{dof}). }

\begin{table}[h!]
  \centering
  \begin{small}
    \begin{tabular}{cccccccccccc}
      \toprule
      Concept $\,$ \textbackslash $\,$ load step & 1-4 & 5-16 & 17-24 & 25-28 & 29-36 & 37-48 & 49-60 \\
      \midrule
      Free \acp{dof} [M] & 0.084 & 0.303 & 0.528 & 0.916 & 1.582 & 2.758 & 4.730  \vspace{0.12cm} \\
      Total background mesh cells [M] & 0.031 & 0.147 & 0.249 & 0.421 & 0.708 & 1.190 & 2.00 \vspace{0.12cm} \\
      Active cells (\% Total cells) & 92.82 & 92.27 & 92.64 & 92.67 & 93.39 & 94.15 & 94.90 \vspace{0.12cm} \\
      Accumulated \ac{amr} cycles & 0 & 3 & 4 & 5 & 6 & 7 & 8  \vspace{0.12cm} \\
      Average nonlinear iters. per load step & 15 & 30 & 12 & 9 & 25 & 9 & 18 \vspace{0.12cm} \\
      \bottomrule
    \end{tabular}
  \end{small}
  \caption{Evolution of \ac{amr} cycles, size of the problem and (average) nonlinear iterations (\VoidStr{} structure).}
  \label{tab:SummaryBullets}
\end{table}

\added{The number of cells and \acp{dof} grow at a rate determined by the parameters of \Alg{ilp} and the growth of the error indicator; see \Fig{IndicatorsError} and \Tab{SummaryBullets}.
Mesh adaptations occur at steps 4, 16, 24, 28, 36 and 48. During the load steps in-between adaptations the mesh remains unchanged. Some relevant conclusions can be extracted out of these results. 
First, these results confirm that, as far as an error indicator-based mesh acceptability criterion is to be fulfilled, there is actually no need to adapt/redistribute the mesh at each load step, as, e.g., it is done in \cite{Frohne2016}. This in turn translates into significantly reduced computational times compared to other heuristic, non-error-driven strategies. This is indeed one of the main goals of \Alg{ilp}. Second, as $\theta_r=0.1$, the ratio of active cells in the meshes corresponding to two consecutive adaptation cycles asymptotically converges to approximately 1.7, as only one mesh adaptation cycle is required to keep $\gre$ below $\mgre$ when \ac{amr} is triggered.
Third, during the whole simulation,
the number of active cells is very close to the number of total cells. Some external cells
are refined to satisfy the 2:1 balance constraint but their number is actually very small. 
Thus, these results confirm that, as far as one coarsens as much as possible exterior cells right at the beginning, they 
do not actually represent a source of overhead for the combined 
use of the unfitted method and \ac{amr}. Finally, in \Tab{SummaryBullets} we also report the average number of nonlinear Newton iterations required to solve the
nonlinear problems at lines \ref{line:solve-nlp-first} and \ref{line:solve-nlp-second} in \Alg{ilp}. This number is
an indication of the difficulty of solving the nonlinear problem on a given mesh and, as it can be seen
from \Tab{SummaryBullets}, quite variable depending, among other factors, on the mesh adaptation patterns and the loading process.}

\subsection{3D Cantilever Beam strong scalability test} \label{sec:CBT-3D}

\begin{figure}[h!]
  \centering
  \includegraphics[width=0.6\textwidth]{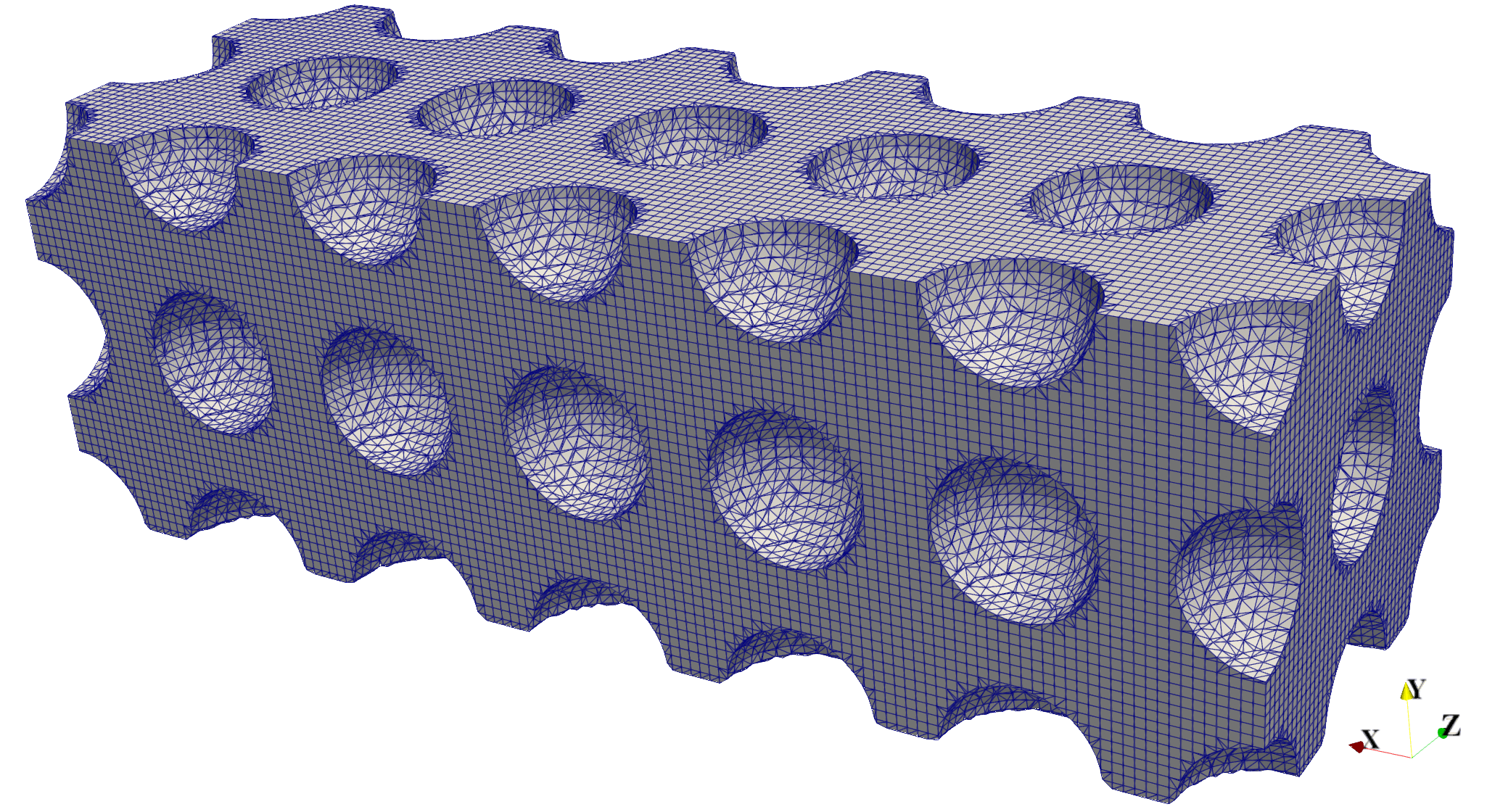}
  \caption{Geometry and initial mesh of the 3D Cantilever beam.}
  \label{fig:CBT-3D-initial-FE-mesh}
\end{figure}

\begin{figure}[h!]
  \centering
  \begin{subfigure}{0.4\textwidth}
    \centering
    \includegraphics[width=\textwidth]{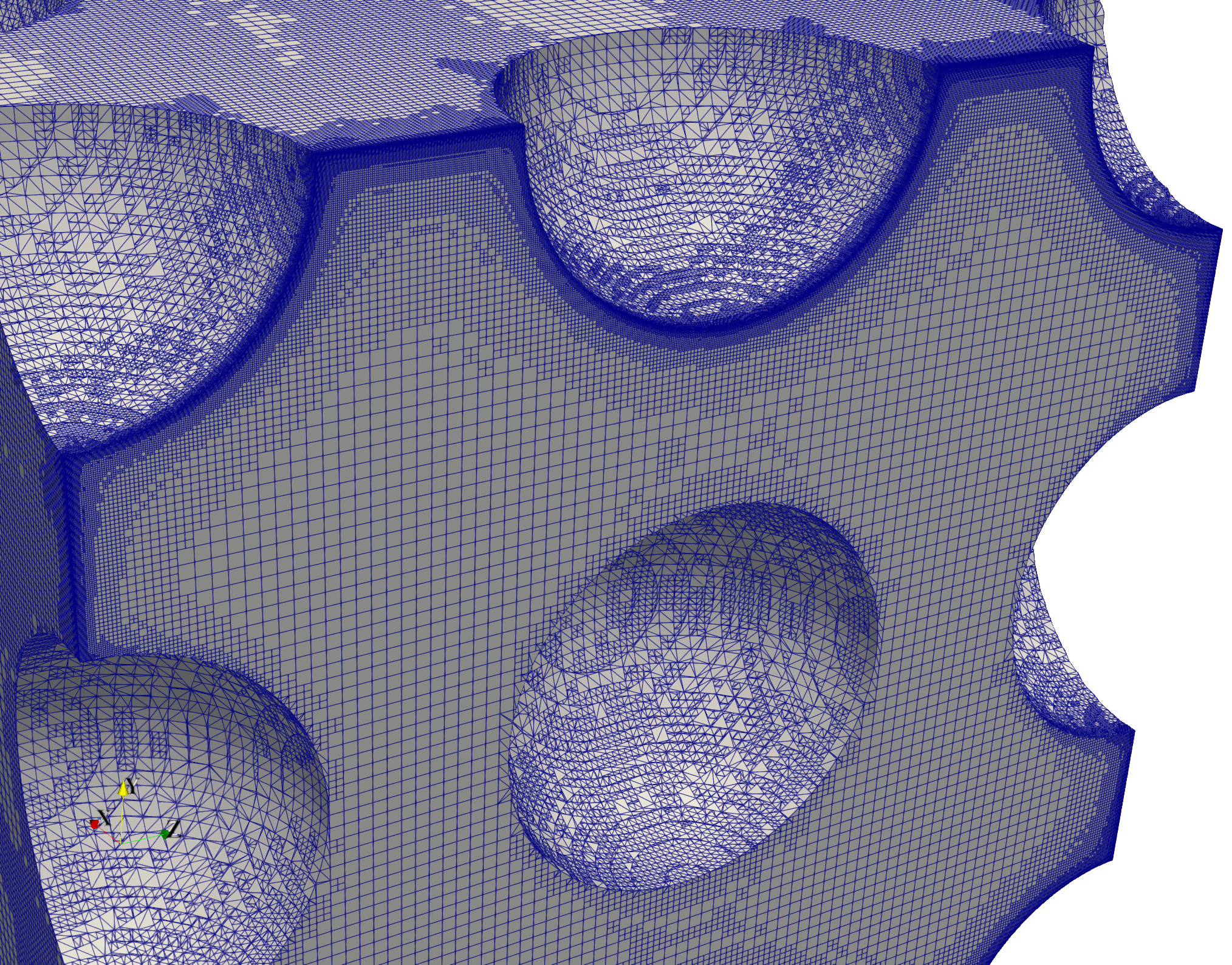}
    \caption{}
    \label{fig:CBT-3D-final-FE-mesh-clamped-face}
  \end{subfigure}
  \begin{subfigure}{0.4\textwidth}
    \centering
    \includegraphics[width=\textwidth]{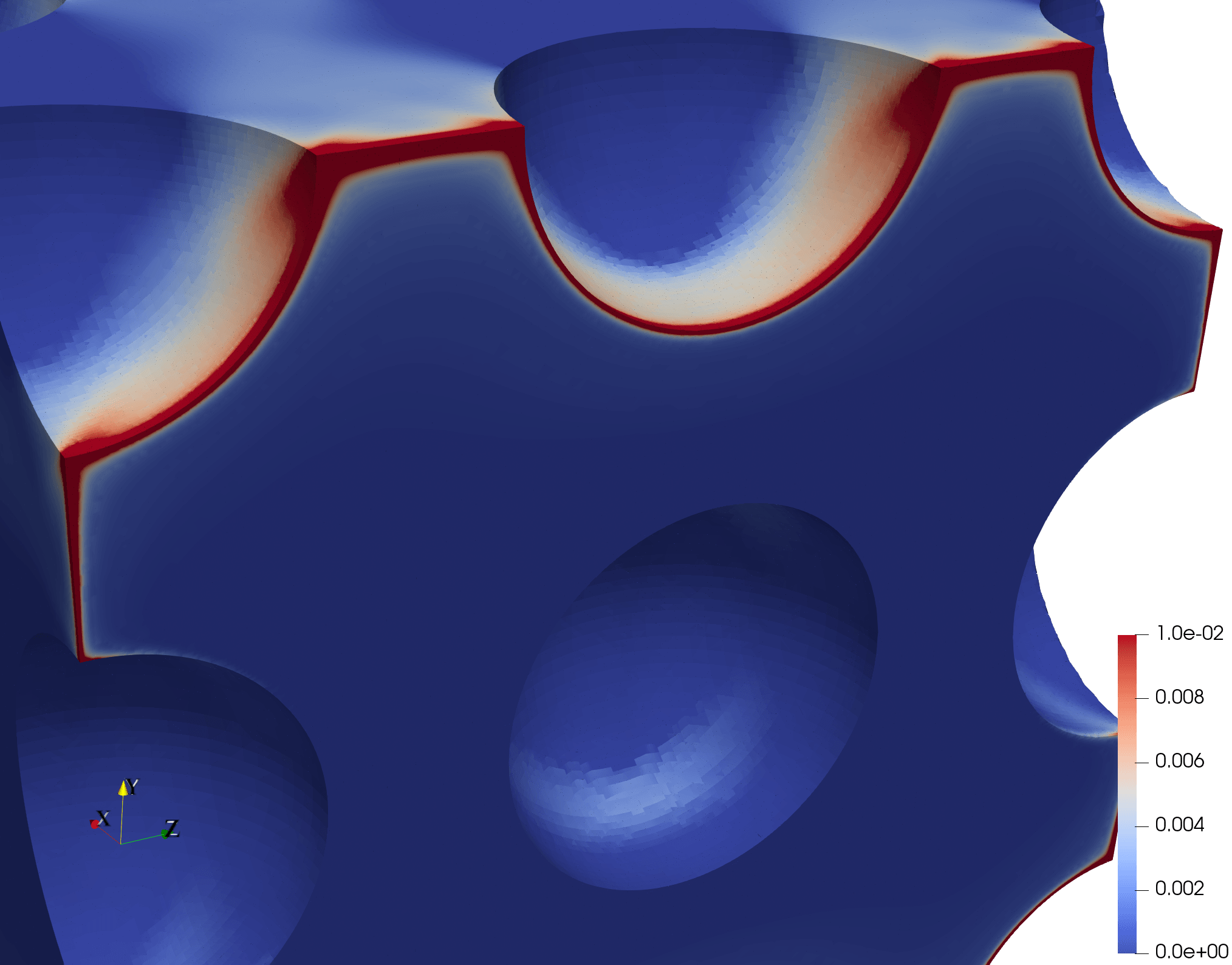}
    \caption{}
    \label{fig:CBT-3D-final-alpha}
  \end{subfigure}
  \caption{3D Cantilever beam test: Final stage: (A) \ac{fe} mesh. (B) J2 Plasticity internal variable $\alpha$.}
  \label{fig:CBT-3D-final-config}
  \end{figure}

\added{
In this section we present the results of a shear test of a geometrically complex
3D beam. The objective of this section is to evaluate the computational performance of \Alg{ilp} for problems with localized inelastic behavior. In particular, we are interested in evaluating and justifying the ability of \Alg{ilp} to efficiently exploit increasing computational resources (cores and memory) when deployed on a parallel supercomputer. 

The problem consists of finding the displacement of a beam subject to a uniform
load applied on the upper surface.  The geometry is a homogeneous medium whose size is $ L \times L/3\times L/3$ with spherical hollows of radius $0.116 \, \text{m}$, arranged in a \ac{sc} lattice of side $L/6$ \cite{Kittel2004}, shown in \Fig{CBT-3D-initial-FE-mesh}, which also shows the initial mesh.  The material properties correspond to a standard A-36 steel having $E=200 \, \text{GPa}$ and $\nu=0.26$. The inelastic behavior is modeled with a linear isotropic J2 plasticity constitutive model, where $K_{inf}=K_0=0$ and $\delta=0$, and the yield threshold and the isotropic hardening factors being $\sigma_{y}=0.250 \, \text{GPa}$ and $H=0.2 \, \text{GPa}$, resp.  Homogeneous Dirichlet boundary  conditions are prescribed in $x=0$ and the rest of the faces are traction free, except the upper one (i.e. $y=L/3$) in which a linearly increasing traction from $0$ to $20 \text{KN}/\text{m}^2$ is applied, discretized into $41$ load steps. We note that the imposition of this Neumann-type boundary condition requires the second integral in the right-hand side of \Eq{residual-force} and thus involves the evaluation of integrals over cut faces. In \Alg{ilp}, we set  $\theta_r=0.07$, $\theta_c=0.03$, $\mgre=0.06$, \ASF=2, and \NAS{}=4. Experimentation with this particular set of parameter values when applied to the same problem but coarser resolution meshes revealed that these were able to achieve a reasonable trade-off among all factors involved (compared to other possible a priori sensible choices). We show the localization at the maximum load step close to the clamped face in \Fig{CBT-3D-final-FE-mesh-clamped-face} (mesh refinement) and \ref{fig:CBT-3D-final-alpha} (J2 isotropic plasticity history variable $\alpha$).}

\begin{figure}[h!]
  \centering
  \includegraphics[width=0.9\textwidth]{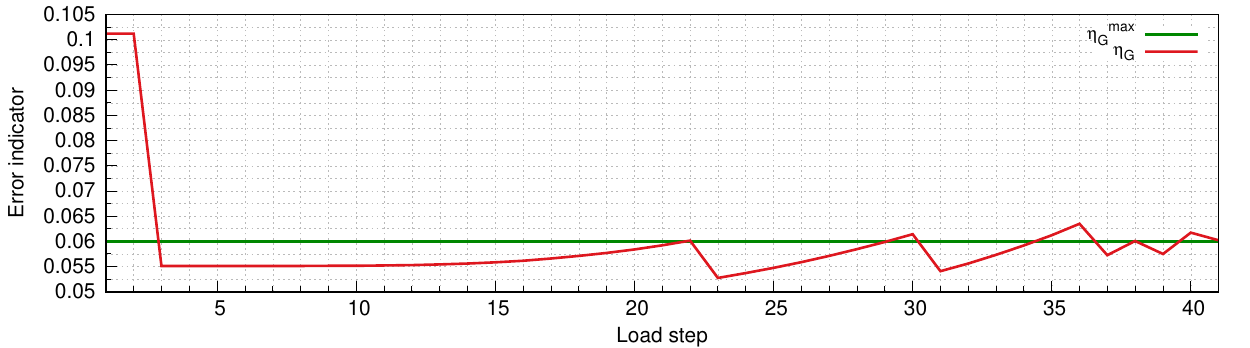}
  \caption{3D Cantilever beam test: Evolution of $\gre$.}  
  \label{fig:CBT-3D-evol-gre}
\end{figure}

\added{
In \Fig{CBT-3D-evol-gre} and \Tab{CBT-3D-summary} we report the same concepts as its counterparts \Fig{IndicatorsError} and \Tab{SummaryBullets}, resp., for the experiment at hand. Overall, very similar conclusions to the ones in the previous section, with some subtle differences, can be extracted from these results.
First, in \Fig{CBT-3D-evol-gre}, we can observe larger intervals in-between  mesh adaptations than those observed in the previous section, and thus even more reduced computational times compared to other heuristic, non-error-driven strategies. Second, in \Tab{CBT-3D-summary} we can also observe a high variability in the number of nonlinear iterations with load step. However, this time, there is a clear steady increase in this number with increasing load.}

\begin{table}[h!]
	\centering
	\begin{small}
		\begin{tabular}{ccccccccccccc}
			\toprule
			Concept $\,$ \textbackslash $\,$ load step & 1-2 & 3-22 & 23-30 & 31-36 & 37-38 & 39-40 & 41 \\
			\midrule
			Free \acp{dof} [M] & 0.592 & 1.637 & 2.425 & 3.553 & 5.258 & 7.840 & 11.69 \vspace{0.12cm} \\
			Total background mesh cells [M] & 0.249 & 0.722 & 1.050 & 1.537 & 2.255 & 3.322 & 4.907 \vspace{0.12cm} \\
			Active cells (\% Total cells) & 82.17 & 90.08 & 91.67 & 92.88 & 94.11 & 94.94 & 95.56 \vspace{0.12cm} \\
			Accumulated adaptivity cycles & 0 & 3 & 4 & 5 & 6 & 7 & 8   \vspace{0.12cm} \\
            Average nonlinear iters. per load step & 1 & 5 & 7 & 10 & 16 & 26 & 33 \vspace{0.12cm} \\
            \bottomrule
		\end{tabular}
	\end{small}
	\caption{Evolution of \ac{amr} cycles, size of the problem and (average) nonlinear iterations (3D Cantilever test).}
	\label{tab:CBT-3D-summary}
\end{table}

\begin{figure}[h!]
  \centering
  \begin{subfigure}{0.45\textwidth}
    \centering
    \includegraphics[width=\textwidth]{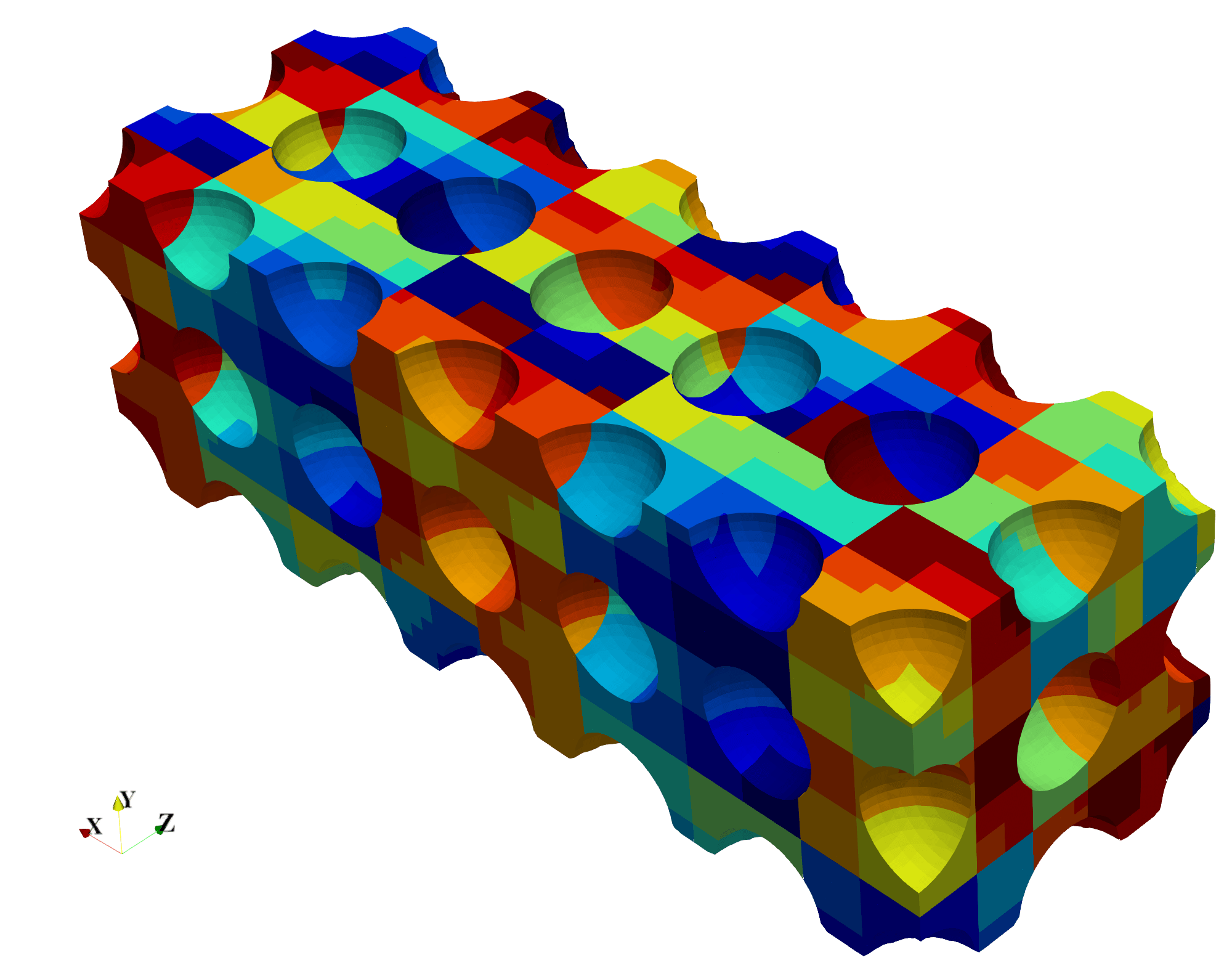}
    \caption{}
    \label{fig:CBT-3D-initial-distribution-of-tasks}
  \end{subfigure}
  \begin{subfigure}{0.45\textwidth}
    \centering
    \includegraphics[width=\textwidth]{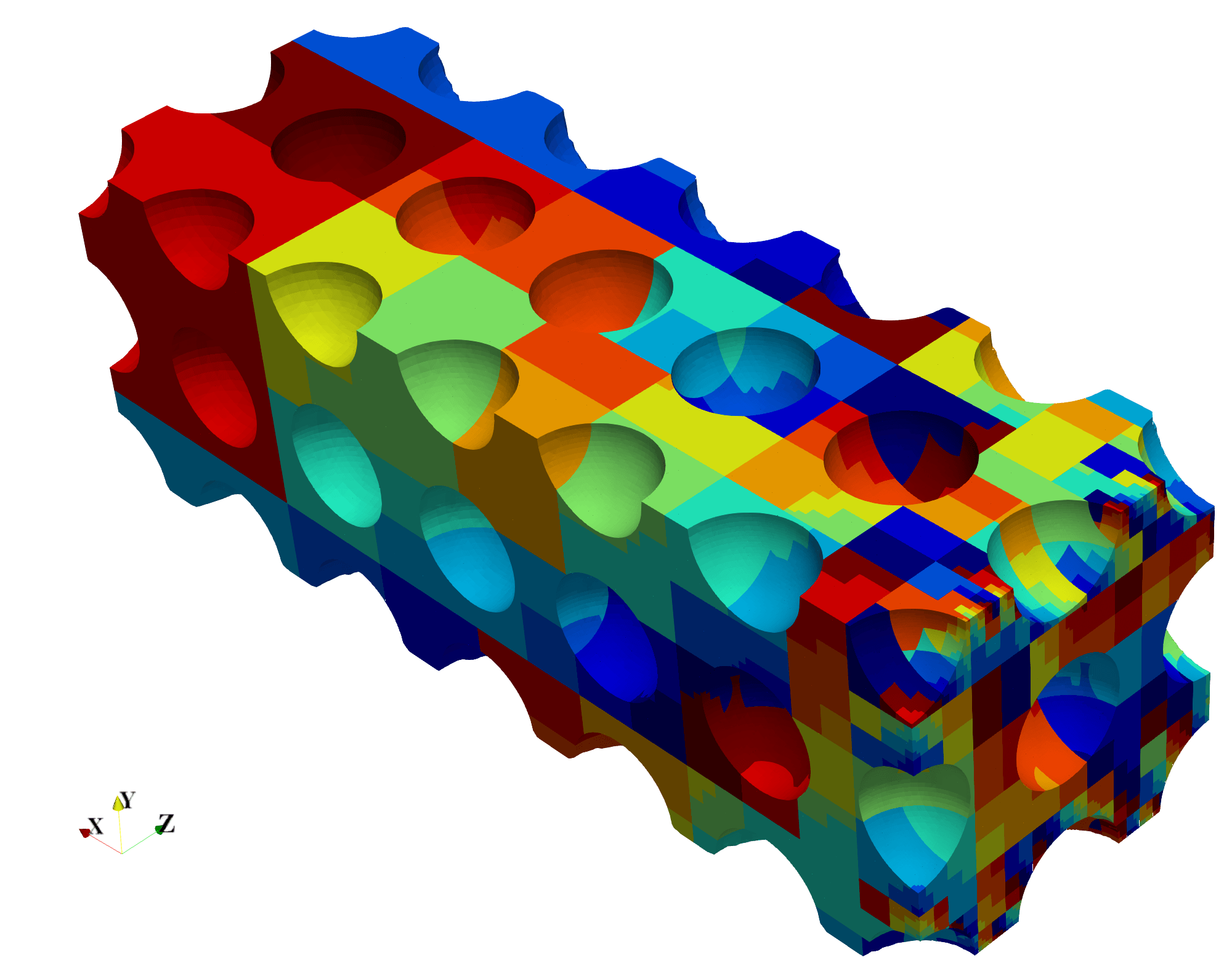}
    \caption{}
    \label{fig:CBT-3D-final-distribution-of-tasks}
  \end{subfigure}
  \caption{3D Cantilever beam: (A) initial and (B) final  distribution of active cells among cores.} 
  \label{fig:CBT-3D-distribution-of-tasks} 
  \end{figure}

\added{We perform a strong scalability study of \Alg{ilp} when applied to the problem described above on the {\em NCI-Gadi}  supercomputer \cite{NCIgadi}; see Section \ref{sec:exp-environment}. The challenge resides in the fact that, as seen in \Fig{CBT-3D-final-config}, localization, and thus refinement along the domain, is not aligned with the initial uniform distribution of load (note that we start the loading process with a uniform mesh). As a consequence, during a load stepping adaptive simulation, the computational load of the different cores will become rapidly highly unbalanced if one does not dynamically reshuffle the mesh and associated fields among them during the simulation. Recall that, in \Alg{ilp}, this is triggered in line~\ref{line:redistr-tasks}. In order to illustrate this, 
in \Fig{CBT-3D-initial-distribution-of-tasks} and \ref{fig:CBT-3D-final-distribution-of-tasks} we show the distribution of load (active cells) among 288 cores in the initial and final load steps. As expected, the initial distribution of load is uniform, with similarly shaped subdomains. However, the final distribution of load is far from uniform, with a subdomains density (number of subdomains per unit volume) much higher close to the clamped face. Indeed, as expected, the cells which cover the $x \geq L/2$ half of the beam are owned by a very small portion of subdomains as effect of dynamic load balancing. We note that, as there is communication overhead involved in the repartitioning of data, it is thus interesting to study whether for the problem at hand this overhead is outweighed by the gains derived from a better distribution of computational load.}

The results of the strong scalability test are given in \Tab{CB-3D-SST-summary}.  
We recall that, as shown in \Tab{CBT-3D-summary}, the size of the discrete problem in the last load step is  around 11.7M~\acp{dof}. The simulation consumes, at its maximum, 
around 86\% of the total amount of memory available when deployed on 4 nodes (192 \added{cores}) on the \textit{NCI-Gadi} supercomputer; see \Sec{exp-environment}. 
 The wall clock times in the table
reported include the full load increment simulation, end-to-end, as reported by the Job Queuing System. Thus, this is the true time that the user experiences when using the software. The value given in the row labeled as EST (\textit{Extrapolated Sequential Time}) is computed as the actual time measured with 192 \added{cores}, multiplied by a factor of 192, and can be considered as an estimation of the sequential execution time under ideal conditions (linear speed-up). In practice, one may expect it to be even higher (as there is typically loss of parallel efficiency to some extent).  Note that we had to estimate this time as the problem at hand does not fit into a single node. 

\begin{table}[h!]
	\centering
	\begin{small}
		\begin{tabular}{rrrrr}
			\toprule
             P & Wall clock Time & $S_p$ & $E_p$ & \#\acp{dof} per \added{core} \\
			\midrule
             EST            &     4936h 32m 00s    &      ---       &      ---       &  11,687,673  \\
             192            &       25h 42m 40s    &      ---       &      ---       &      60,873  \\
             384            &       10h 46m 22s    &     2.390      &     1.190      &      30,436  \\
             768            &        5h 24m 15s    &     4.760      &     1.190      &      15,218  \\
             1,536           &        2h 58m 15s    &     8.650      &     1.080      &       7,609  \\
             3,072           &        1h 49m 14s    &    14.120      &     0.880      &       3,804  \\
            \bottomrule
		\end{tabular}
	\end{small}
    \caption{3D Cantilever beam - strong scalability test. \textit{\#\acp{dof} per \added{core}} is computed as the total number of \acp{dof} divided by the number of \added{cores} ($p$). Simulation ends after 41 load steps.}
	\label{tab:CB-3D-SST-summary}
\end{table}

The column labeled as $S_p$ (speed-up) in \Tab{CB-3D-SST-summary} is defined as the ratio between the parallel execution time on $192$ \added{cores} (reference time), referred to as $t_{192}$, and the parallel execution time on $p$ \added{cores}, i.e. $S_p=t_{192}/t_p$. The column labeled as $E_p$ (parallel efficiency), is defined as 
$E_P=S_p/k$, with $k=p/192$. Clearly, 
the most salient property of the framework, as can be shown in 
\Tab{CB-3D-SST-summary}, is its ability to efficiently reduce execution times with increasing number of cores despite \added{the localization challenge.
A closer inspection of the load balance achieved reveals that the 
relative difference among the subdomains with the maximum and minimum number of \acp{dof} stays well below 1\% during the whole simulation versus differences of more than 100\% when dynamic load balance is deactivated.} By exploiting parallel resources, the execution time is reduced by a factor of 
$14.120$ when $3072$ \added{cores} are used. For practical purposes, these speed-ups reflect a reduction of the execution time in 
industrial simulations from months (as in the EST) to 
less than a couple of hours when $3$K \added{cores} are used. The loss of parallel performance with $p$ is associated to parallelism related overheads; in this case, more computationally intensive simulations, involving larger
loads per core, are required to exploit the computational resources more 
efficiently.

In order to gain more insight on the results in \Tab{CB-3D-SST-summary}, we profiled the execution of the framework for $P=192$ \added{cores}. It turns out that the bulk of the computation is concentrated in the nonlinear solver, which amounts to 75\% of the total computation time (i.e. aggregated across all load steps). Out of this 75\%, the computation of the residual, Jacobian, and preconditioned iterative solution of linear systems concentrate 32\%, 14\%, 28\%, resp., of the total computation time. This is followed by the computation of error estimators, which concentrates 9\% of the total computation time. Finally,
the updates of global \ac{fe} spaces required after adaptation and redistribution
amounts to 3\%, and those related to the background mesh, embedded domain boundary intersection, and setting up the aggregates, to 0.5\%. The timings in \Tab{CB-3D-SST-summary}, both in magnitude, and the rate at which they decrease with $P$ (note that the computation of the residual and Jacobian are highly parallel stages), can be justified, among others, by the number of nonlinear iterations required to achieve convergence; see results in Table~\ref{tab:CBT-3D-summary}.
\added{Besides, the acceleration caused by sticking into a balanced computational load largely outweighs the extra overhead associated to dynamic load balancing.}

\section{Conclusions}\label{sec:conclusions}

This work extends the \ac{hagfem} to the 
context of nonlinear solid mechanics problems.
To the best of the authors' knowledge, this is the first fully parallel  distributed-memory unfitted $h$-adaptive \ac{fe} framework robust with respect to small cut
cells that solves nonlinear solid mechanics problems. It is grounded on five 
main building blocks:
1) unfitted formulations (using aggregation techniques) that are robust irrespective of the cut locations, 
2) a strategy to deal with history variables of any tensorial order provided by the constitutive 
model in combination with aggregated spaces,  
3) an algorithm for the solution of nonlinear problems based on the Newton-Raphson method, endowed with optimization strategies to improve
the performance of the nonlinear solver in highly nonlinear scenarios, 
4) an algorithm to solve the nonlinear problem incrementally that includes hierarchical \ac{amr} to provide an adapted mesh able to 
capture the evolution of inelastic behavior and 
5) a distributed-memory implementation that relies on parallel linear solvers and octree engines with dynamic 
load-balancing at each adaptive step.

The proposed framework has been 
extensively tested with a large set of experiments
for both irreducible and \mixedup{} formulations that include incompressible materials. 
These experiments reveal good accuracy in terms of error 
convergence when comparing the numerical solution with the analytical one (when provided) 
and show that the onset of nonlinearity in 
unfitted boundaries is properly captured; this becomes a crucial issue in tests 
involving complex geometries where the inelastic behavior takes place in weak 
regions of the physical domain. These
experiments not only have verified 
efficiency and robustness of the framework, but also revealed the sensitivity of 
the \ac{amr} strategy to capture localized nonlinear regions as consequence of 
complex load scenarios. The scalability properties of the framework 
have also been analyzed.

We have applied this framework to more complex periodic structures in order 
to show how the framework can efficiently deal with these geometries, avoiding the need to generate very costly body-fitted meshes and graph partitioners. Body-fitted meshing is a computational bottleneck that requires human intervention and prevents an automatic geometry-to-solution workflow. On the contrary, the framework proposed in this work is fully automatic and highly scalable.
This fact enables, e.g. robust optimization of structures
or uncertainty quantification of structures with random
domains \cite{BadiaHampton2021}. This framework is particularly relevant in
\ac{am}, in order to virtually certify and optimize mesoscale lattice structures
or to quantify the effect of geometrical irregularities, or to handle growing geometries in \ac{am} process simulations \cite[]{Neiva2019}.

\section*{Acknowledgments}

\thethanks

\begin{small}

\setlength{\bibsep}{0.0ex plus 0.00ex}
\bibliographystyle{myabbrvnat}
\bibliography{art044.bib}

\end{small}

\end{document}